\documentclass[a4paper,11pt]{article}
\usepackage{amssymb}
\usepackage{cite}
\usepackage{color}

\newcommand{\s}{s}

\newcommand{\cL}{{\cal L}}

\newcommand{\cC}{\mathcal{C}}
\newcommand{\cB}{\mathcal{B}}
\newcommand{\cD}{\mathcal{D}}

\newcommand{\cK}{\mathcal{K}}

\newcommand{\R}{\mathbb{R}}

\newcommand{\N}{\mathbb{N}}

\newcommand{\cF}{{\cal F}}
\newcommand{\cuad}{{\sqcap\kern-.68em\sqcup}}

\newcommand{\norm}[1]{\|#1\|}
\newcommand{\equ}[1]{(\ref{#1})}

\newtheorem{theorem}{Theorem}[section]
\newtheorem{proposition}{Proposition}[section]

\newtheorem{lemma}{Lemma}[section]
\newtheorem{corollary}{Corollary}[section]

\begin{document}

\begin{center}{\bf     On  isolated singularities for fractional Lane-Emden equation    \\[1.5mm]

  in the Serrin's critical case

 }\bigskip

 {\small \begin{center}

{\small
 Huyuan Chen \footnote{chenhuyuan@yeah.net}  \quad{and}\quad  Hichem Hajaiej \footnote{hhajaie@calstatela.edu}

\medskip
$^1$Department of Mathematics, Jiangxi Normal University,\\
Nanchang, Jiangxi 330022, PR China \\[10pt]

$^2$Department of Mathematics, College of  Natural Science
Cal State University,\\ 5151 State Drive, 90032 Los Angeles, California, USA  \\[19pt]
 }

\end{center}
}

\begin{abstract}

A challenging problem for nonlocal semilinear  elliptic equation is the fractional Lane-Emden equation in the Serrin's critical case. This is due to the various challenges and the lack of tools to analyze the isolated singularities. These difficulties come from the nonlocal setting where the ODE's tools fail to provide useful information. Additionally, the computation of the fractional Laplacian for the involved logarithmic functions is tricky. Finally, it is not an easy task to show that the solutions are in the appropriate function space. In this paper, we solve the fractional Lane-Emden equation in the Serrin's critical case for the fractional Laplacian by developing an innovative and self-contained approach that also applies to the classical setting ( Laplacian). 

We give a classification of the isolated singularities of positive  solutions to the semilinear fractional elliptic  equations
$$(E) \qquad\qquad (-\Delta)^s u =  u^{\frac{N}{N-2s}}\quad {\rm in}\ \ \Omega\setminus\{0\},\qquad u\geq 0\quad{\rm in}\ \ \R^N\setminus\Omega,\qquad\qquad\quad   $$
where $s\in(0,1)$,  $\Omega$ is a bounded  domain containing the origin in $\R^N$ with $N>2s$ and $\frac{N}{N-2s}$ is the Serrin's critical exponent. We use an initial asymptotic at infinity to transform the critical case into a subcritical case where the underlying equation involves the fractional Hardy operator. The construction of singular solutions is based on the fact that some special functions are subsolutions of the original problem near the origin.

We also classify the non-removable singularities of the solutions of $(E)$ and show the existence of a sequence of  isolated singular solutions parameterizing the coefficients of the second order blow up term.  To the best of our knowledge, this idea has also been first used in this work.
\end{abstract}

  \end{center}
  % \tableofcontents \vspace{1mm}
   \smallskip

   \noindent {\small {\bf AMS Subject Classifications}:     35R11;  35J75; 35A01.
  \smallskip

 %\noindent {\small {\bf Acknowledgment}: The authors would like to thank NYU Shanghai where they spent two years. Special thanks to Pr. Fanghua Lin for fruitful and inspiring discussions. This work has been initiated during the first year of their visit.

\vspace{2mm}

\setcounter{equation}{0}
\section{Introduction}
   
   Let $s\in (0,1)$, $\Omega$ be a bounded $C^2$ domain containing the origin in $\R^N$ with $N$ a positive integer such that $N>2s$ and
 $$p^*=\frac{N}{N-2s}$$
  is the Serrin's critical exponent.
  The purpose of this paper is to classify  the isolated singular  positive solutions of the elliptic problem
\begin{equation}\label{eq 1.1}
\left\{
\begin{array}{lll}
(-\Delta)^s u=   u^{p^*}  \quad &{\rm in}\ \, \Omega \setminus\{0\},
\\[2mm]
\qquad \ \ u=h&{\rm in}\ \,  \R^N\setminus \Omega,
\end{array}\right.
 \end{equation}
where  $(-\Delta)^s$ is the fractional Laplacian   defined by
$$
(-\Delta)^\s  u(x)= c_{N,\s}\lim_{\epsilon\to0^+} \int_{\R^N\setminus B_\epsilon }\frac{ u(x)-
u(x+z)}{|z|^{N+2\s}}  dz,
$$
   $B_r(y)\subset \R^N$ is the ball of radius $r>0$ centered at $y$. Here and in what follows,   $B_r=B_r(0)$,
   $$c_{N,\s}=2^{2\s}\pi^{-\frac N2}\s\frac{\Gamma(\frac{N+2\s}2)}{\Gamma(1-\s)}$$
   is  the normalized constant, see \cite{NPV},  with $\Gamma$ being the Gamma function.  It is known that  $(-\Delta)^\s  u(x)$ is well defined if $u$ is twice continuously differentiable in a neighborhood of $x$, and contained in the space
$L^1_s(\R^N):=L^1(\R^N,\frac{dx}{1+|x|^{N+2\s}})$.    We also note that for $u \in L^1_s(\R^N)$ the fractional Laplacian $(-\Delta)^\s  u$ can also be defined as a distribution:
$$%\begin{equation}
 % \label{eq:-distributional-sense}
 \langle(-\Delta)^\s u, \varphi \rangle= \int_{\R^N} u (-\Delta)^\s \varphi \,dx \quad\  {\rm for\ any }\ \  \varphi \in C^\infty_c(\R^N).
$$%\end{equation}
We then have
$$%\begin{equation}
 % \label{eq:-distributional-sense}
\cF((-\Delta)^\s u) = |\cdot|^{2\s} \hat u \qquad  {\rm in}\ \  \R^N \setminus \{0\}
$$%\end{equation}
in the sense of distributions,  both $\cF(\cdot)$ and $\hat \cdot$ denote the Fourier transform. \smallskip

When $s=1$ and $N\geq 3$, the isolated singularities of the solutions of the related model equation
\begin{equation}\label{eq 1.1-loc}
\left\{
\begin{array}{lll}
 -\Delta  u=   u^{p}  \quad &{\rm in}\ \, \Omega \setminus\{0\},
\\[2mm]
\quad \ \ u=0&{\rm on}\ \,  \partial \Omega
\end{array}\right.
 \end{equation}
 were studied extensively in the last decades.  When $p\in(1,\, \frac{N}{N-2})$, Lions in \cite{L}
 classified the singular solutions by establishing a connection with the weak solutions of
 $$-\Delta  u=   u^{p} +k\delta_0 .   $$
 He showed that the  positive  solutions
 have  the asymptotic behavior
 $$\lim_{|x|\to0^+} u_k(x)|x|^{N-2}=c_Nk $$
 for some $k\geq0$, and a normalized constant $c_N>0$. In the particular case $k=0$, the solutions have removable singularity at the origin.
This is always the case when $p\geq \frac{N}{N-2}$. Therefore this method fails to classify the singularities in this case.   When $p\in\big(\frac{N}{N-2},\, \frac{N+2}{N-2}\big)$, the singularity of the solutions to (\ref{eq 1.1-loc}) was studied by Gidas-Spruck in \cite{GS} by using analytic techniques from \cite{A2}. In this case, the positive singular solutions have the following type of singularity:
 $$u(x)= c_p |x|^{-\frac{2}{p-1}}(1+o(1))\quad{\rm as}\ \, |x|\to0^+$$
 with  the coefficient $c_p^{p-1}=\frac{2}{p-1}(N-2-\frac{2}{p-1}).$
 When $p=\frac{N}{N-2}$, the author in \cite{A1,A2} gave the following classification:  Any positive solution $u$ of (\ref{eq 1.1-loc}) has a removable singularity at the origin or has the asymptotic behavior
 $$u(x)=\cK_1 |x|^{2-N}(-\ln|x|)^{\frac{2-N}{2}}(1+o(1))\quad{\rm as}\ \, |x|\to0^+,$$
where
$$\cK_1=\big(\frac{(N-2)^2}{2}\big)^{\frac{N-2}{2}}.$$
Moreover,   the existence of a singular solution is obtained by the Phase plane analysis.  Using these solutions, Pacard in\cite{P} constructed positive solutions with the prescribed singular set.
 Later on, Caffarelli-Gidas-Spruck in \cite{CGS} (also see \cite{KM}) classified the singular solutions of  (\ref{eq 1.1-loc}) when $p=\frac{N+2}{N-2}$. In this case, they satisfy
 $$u(x)=\varphi_D(\ln |x|) |x|^{-\frac{N-2}{2}}(1+o(|x|))\quad{\rm as}\ \, |x|\to0^+,$$
  $\varphi_D$ may be a constant $c_{\frac{N+2}{N-2}}$ or a periodic function. We refer to \cite{BV,HLT,S0,S} for more  details for elliptic problem with general differential operators.

 Due to the numerous applications of Dirichlet problems with nonlocal operators,  there has been an increasing interest in this topic during the last years.  A general setting of nonlinear nonlocal operators  are studied in \cite{CS0,CS1,CS2}, where the authors connect the fractional problem with an underlying second order  degenerated problem in a half space in a higher dimension. For regularity properties of the solutions of the fractional problem, the reader could see \cite{CS2,RS,Ls07}. For  blow-up analysis for nonlocal problems, see  \cite{CJ,CV1,CFQ}.
 
Motivated by \cite{L},  the isolated singularity of
 \begin{equation}\label{eq 1.2}
\left\{
\begin{array}{lll}
(-\Delta)^s u=   u^{p}  \quad &{\rm in}\ \, \Omega \setminus\{0\},
\\[2mm]
\qquad \ \ u=0&{\rm in}\ \,  \R^N\setminus \Omega
\end{array}\right.
 \end{equation}
was studied in \cite{CQ} for $p\in(1, p^*)$ via the connection with the distributional solutions  of
$$ \left\{
\begin{array}{lll}
(-\Delta)^s u=   u^{p} +k\delta_0  \quad &{\rm in}\ \, \Omega \setminus\{0\},
\\[2mm]
\quad \ \ u=0&{\rm in}\ \,  \R^N\setminus \Omega,
\end{array}\right.
 $$
where the positive solution with singularity can be described by fundamental solution of fractional Laplacian.
 As in the Laplacian case, this method does not apply to classify isolated singularities for the Lane-Emden equation with Serrin's critical and super critical exponent, i.e. $p\geq p^*$.
 When $p\in(p^*,\frac{N+2s}{N-2s}]$,  \cite{JLX,CJ} built the platform of the isolated singularity for positive solution to
  \begin{equation}\label{eq 1.3}
\left\{
\begin{array}{lll}
 \,\quad\ \  {\rm div}(t^{1-2s}) \nabla U)=0  \quad &{\rm in}\ \, B_1\times (0,1)
\\[2mm]
\displaystyle  \lim_{t\to0^+}t^{1-2s}\partial_t U(x,t) =U^{p}(x,0)&  {\rm on}\ \,   B_1\setminus\{0\},
\end{array}\right.
 \end{equation}
 in the Sobolev critical case $p=\frac{N+2s}{N-2s}$. The non-removable singular solutions of (\ref{eq 1.3}) behave as follows:
 $$c_2 |x|^{-\frac{N-2s}{2}}\leq u(x) \leq c_1|x|^{-\frac{N-2s}{2}}$$
  for some $c_1>c_2>0$, \cite{YZ} gives a   description of singular solutions for $p\in(p^*,\frac{N+2s}{N-2s})$
 $$c_4|x|^{-\frac{2s}{p-1}}\leq u(x)\leq c_3|x|^{-\frac{2s}{p-1}}$$
 for some $c_3>c_4>0$.   For the Serrin case $p=p^*$ \cite{WW} shows the bound 
   \begin{eqnarray*}
b_1(|x|(-\ln|x|)^\frac{1}{2s})^{2s-N} \leq  u(x) 
  \leq  b_2(|x|(-\ln|x|)^\frac{1}{2s})^{2s-N} 
\end{eqnarray*}
for some $0<b_1<b_2<+\infty$.

   It is worth noting that (\ref{eq 1.3}) is a local degenerated problem.  For the existence of isolated singular solutions,  \cite{AC} constructed a sequence of isolated singular solutions of (\ref{eq 1.2}) with $\Omega=\R^N$ and $p\in(p^*,\frac{N+2s}{N-2s})$, with fast decaying at infinity, (also see \cite{AD,DG} for  the fractional Yamabe problem with isolated singularities in the Sobolev critical case).

It is known that   the isolated singularity of solutions to (\ref{eq 1.2}) in the  Serrin's critical case is totally different, since it doesn't blow up with a negative-power function, therefore  it is much more difficult to analyze the isolated singularities.  The main challenges arise from the nonlocal property: One is that the solutions must remain in $L^1_s(\R^N)$ in some type of scaling in the blow-up analysis; the second  is the tricky calculation of the fractional Laplacian for the involved logarithmic functions; the last and the most challenging is  to obtain the singular solutions for the fractional problem, while we don't have  phrase plane analysis for the existence.

Our aim in this paper is two-fold: the first is to classify the positive singular solutions of the fractional  problem (\ref{eq 1.2}) in the Serrin's critical case and the second is to show the existence of singular solutions.  Let us first introduce some notations.  For the critical case,  we note that   $-\frac{2s}{p^*-1} =2s-N$.   For $\tau\in(-N,2s)$, we denote
 \begin{equation}\label{cons 0}
  \cC_{s}(\tau) =2^{2\s}\frac{\Gamma(\frac{N+\tau}{2})\Gamma(\frac{2\s-\tau}{2})}{\Gamma(-\frac{\tau}{2})\Gamma(\frac{N-2\s+\tau}{2})}
  \end{equation}
 and
 \begin{equation}\label{cons 1}
   \cC_{s}(0)=0,\qquad \cC_{s}'(0)=-2^{2\s-1}\frac{\Gamma(\frac{N}{2})\Gamma(\s)}{ \Gamma(\frac{N-2\s}{2})}  <0.
 \end{equation}

 Let
 \begin{equation}\label{cons k}
  \cK_s=\Big(- \cC_{s}'(0)\frac{N-2s}{2s}\Big)^{\frac{N-2s}{2s}},
 \end{equation}
 then
 \begin{equation}\label{cons k-l}
 \lim_{s\to1^-}\cK_s=\Big(\frac{(N-2)^2}{2}\Big)^{\frac{N-2}{2}}=\cK_1.
 \end{equation}
 See {\it Appendix A} for the details of the proof.\smallskip

The classification of isolated singularities of (\ref{eq 1.1}) states as follows:

\begin{theorem}
\label{teo 5.1}

  Let   $B_1\subset \Omega\subset B_{R_0}$ for $R_0\geq1$,  $h$ be a nonnegative function in $C^\theta(B_{2R_0})\cap L^1_s(\R^N)$ with $\theta>2s$,   
 and   $u$ be a  positive solution of (\ref{eq 1.1}), then $u$ has either a removable singularity at the origin or
  \begin{eqnarray*}
\frac{\cK_s}{C_0} \leq  \liminf_{|x|\to0^+}u(x) (|x|(-\ln|x|)^\frac{1}{2s})^{N-2s}
 &\leq& \cK_s  \\[1.5mm]&\leq&\limsup_{|x|\to0^+}u(x) (|x|(-\ln|x|)^\frac{1}{2s})^{N-2s}\leq  C_0\cK_s,
\end{eqnarray*}
where $C_0\geq 1$ is the best constant in Harnack inequality.
 \end{theorem}

Unlike the Serrin's super critical case, the bound $\displaystyle \lim_{|x|\to0^+}\big( u(x)|x|^{N-2s}\big)<+\infty$  derived by the blow-up analysis is not sharp. In order to improve the upper bound, we use the Liouville property of the fractional Poisson problem
$$\left\{
\begin{array}{lll}
(-\Delta)^s u= \frac{\nu}{|x|^{2s}(-\ln|x|)} u+   f \quad &{\rm in}\ \, B_{R_0} \setminus\{0\},
\\[2mm]
\qquad \ \   u\geq 0&{\rm in}\ \,  \R^N\setminus B_{R_0},
\end{array}\right.
$$
 where $\nu>0$.  For the lower bound, we first obtain a rough one: $\displaystyle\lim_{|x|\to0^+} \big(u(x)|x|^{-\tau}\big)=+\infty$. We then use some special auxiliary tools to improve it to: $$\displaystyle\liminf_{|x|\to0^+} \left(u(x)\big(|x|(-\ln|x|)^\frac{1}{2s}\big)^{N-2s}\right)>0.$$ The improvement of the bound is performed by considering $(-\Delta)^sw$ near the origin, where $w(x)=(|x|(-\ln|x|)^\frac{1}{2s})^{2s-N}$ for $|x|>0$ small, see the details in Proposition \ref{pr 2.1-a} below.

    Because of the nonlocal property,   we can't transform our problem into an ODE and it is a challenging   problem is to get the precise blow-up behavior:
$$\lim_{|x|\to0^+}u(x) \big(|x|(-\ln|x|)^\frac{1}{2s}\big)^{N-2s}=\cK_s. $$
To achieve this goal,   develop new techniques have to be involved. \smallskip

 Our second purpose is to construct isolated singular solutions of (\ref{eq 1.1}).

\begin{theorem}\label{teo 2}
 Let $h=0$ and $\Omega=B_1$,   then there exists $k^*\in\R$ such that  for any $k\in(-\infty,k^*)$
  problem (\ref{eq 1.1}) has  a positive singular solution $u_k$  such that
  $$ \lim_{|x|\to0^+}\Big(u_k(x) -\cK_s \big(|x|(-\ln|x|)^\frac{1}{2s}\big)^{2s-N}\Big)|x|^{N-2s}(-\ln|x|)^{\frac{N}{2s}}=k.    $$
%Moreover,
 %$$\lim_{k\to-\infty} \|u_k\|_{L^1(B_1)}=0 $$
 \end{theorem}

Notice that  a positive solution of (\ref{eq 1.1}) with a removable singularity could be derived by Mountain Pass theorem. Our construction of singular solutions is based on the observation that
some special functions could be sub-solutions  near the origin and   the approximation procedure
is applied to obtain singular solutions by suitable estimates and scaling property.

The derivation of a sequence of  isolated singular solutions   parameterized in the second
 blow-up rate is totally new and our construction is also appropriate to the laplacian case:
\begin{corollary}\label{cr 2}
  Let $p=\frac{N}{N-2}$,    then there exists $k^*\in\R$ such that  for any $k\in(-\infty,,k^*)$
  problem (\ref{eq 1.1-loc}) has  a positive singular solution $u_k$  such that
  $$ \lim_{|x|\to0^+}\Big(u_k(x) -\cK_1 (|x|(-\ln|x|)^\frac{1}{2})^{2-N}\Big)|x|^{N-2}(-\ln|x|)^{\frac{N}{2}}=k.    $$
 
 \end{corollary}

It is important to mention that  the positive singular solutions in general bounded domain can be
derived by approximation and scaling techniques.  On the contrary, in the whole space
problem (\ref{eq 1.1}) has no positive solution. More precisely:

  \begin{proposition}\label{teo 3}
Problem
   \begin{equation}\label{eq 1.1-ext}
(-\Delta)^s u=   u^{p^*}  \quad  {\rm in}\ \, \R^N\setminus\{0\}
 \end{equation}
    has no positive solutions.
\end{proposition}

 The  nonexistence in Proposition \ref{teo 3} isn't new,  \cite[Theorem 1.3]{FQ}.   However,  our methods for all the other results are new and self-contained. In fact, for the Serrin's critical case,  we use some initial asymptotic at infinity to transform the critical problem into a subcritical case for an underlying problem with  the fractional Hardy operator. \smallskip

 The remainder of this paper is organized as follows. Section 2 includes basic tools for the fractional Laplacian, the fundamental calculations of fractional Laplacian on functions involving the logarithmic functions.     In Section 3, we obtain the blow-up upper and lower bounds
 and we improve the isolated singularity to prove Theorem \ref{teo 5.1} in Section 4.
  Section 5 is devoted to the existence of isolated singular solutions and prove Theorem \ref{teo 2}.   Finally, we discuss the constants $\cC_s(0)$ and $\cK_s$,  the  radially symmetric for singular solutions and the nonexistence of positive solution of (\ref{eq 1.1}) in $\R^N\setminus\{0\}$.

\setcounter{equation}{0}
\section{ Preliminary  }

\subsection{Poisson problem }

We start this section by recalling some important facts about the fractional Poisson problem
 \begin{equation}\label{eq 2.1fk}
\left\{
\begin{array}{lll}
(-\Delta)^s u=   f  \quad &{\rm in}\ \, \Omega \setminus\{0\},
\\[2mm]
\qquad \ \  u=0&{\rm in}\ \,  \R^N\setminus \Omega.
\end{array}\right.
\end{equation}
Here and in what follows,   for the sake of simplicity, we  always set
 $$B_1\subset \Omega\subset B_{R_0} $$
 for some $R_0\geq 1$.

\begin{theorem} \label{theorem-C}\cite[Theorem 1.3]{CW}
Let    $f\in C^\theta_{loc}(\bar \Omega\setminus\{0\})$ for some $\theta \in (0,1)$.
 \begin{enumerate}
  \item[(i)] If $f\in L^1(\Omega)$, then for every $k \in \R$ there exists a unique solution
$u_k \in L^\infty_{loc}(\R^N \setminus \{0\}) \cap L^1_s(\R^N)$ of the problem (\ref{eq 2.1fk})
  satisfying the distributional identity
\begin{equation}
  \label{eq:distributional-k-special-case}
\int_{\Omega}u_k   (-\Delta)^s \xi \,dx = \int_{\Omega}f \xi\, dx+ k \xi(0)\qquad   {\rm for\ all}\   \xi \in  C^2_0(\Omega).
\end{equation}
If moreover $f\in L^\infty(\Omega, |x|^\rho dx)$ for some $\rho<2s$, then $u_k$ has the asymptotics
 \begin{equation}\label{beh-1-special-case}
 \lim_{x \to 0} u(x) |x|^{N-2 s} = \frac{k}{c_{ s,0}} \qquad  {\rm with}\quad c_{ s,0}:=c_{N,s} \omega_{_{N-1}} \int^1_0 \int_{B_{t}(0)} \frac{|z|^{2 s-N}-1}{|e_1-z|^{N+2 s}}  dzdt.
 \end{equation}
\item[(ii)] Assume that $f$ is nonnegative and $\int_{\Omega} f\, dx  =+\infty$. Then problem (\ref{eq 2.1fk}) has no nonnegative distributional solution $u \in L^\infty_{loc}(\R^N \setminus \{0\}) \cap L^1_s(\R^N)$.
\end{enumerate}
\end{theorem}

  Motivated by above theorem,  we  can obtain the nonexistence of  solutions for
  \begin{equation}\label{eq 2.2-0}
  \left\{
\begin{array}{lll}
(-\Delta)^s u=   f  \quad &{\rm in}\ \, \Omega \setminus\{0\},
\\[2mm]
\qquad \ \  u=h&{\rm in}\ \,  \R^N\setminus \Omega,
\end{array}\right.
\end{equation}
where  $h\in C^{2}(B_{2R_0})\cap L^1_s(\R^N)$ is nonnegative.

 \begin{theorem}\label{teo 2.1}  Let $f$ be a nonnegative function in $C^\beta_{loc}(\bar \Omega\setminus\{0\})$ with $\beta\in(0,1)$, and $h$ be a nonnegative function  in  $C^2(B_{2R_0})\cap L^1_s(\R^N)$.  

Then  problem (\ref{eq 2.2-0})   has no positive  solution,
  if
 $$\lim_{r\to0^+}\int_{B_{r_0}(0)\setminus B_r(0)}f(x)  dx=+\infty.$$
Particularly, the above assumption could be replaced by
$$\liminf_{|x|\to0^+} f(x)|x|^{N}(-\ln|x|)^{\nu} >0\quad {\rm for\ }\ \nu\leq 1.$$

\end{theorem}

To prove Theorem \ref{teo 2.1}, we need the following comparison principle.

\begin{lemma}\label{cr hp}
Assume that $O$ is a  bounded, Lipschitz continuous domain containing the origin,  $f_i\in C^\beta_{loc}(\bar O\setminus\{0\})$, $h_i\in C^{2}(B_{2R_0})\cap L^1_s(\R^N)$ and $u_i$ with $i=1,2$ are classical solutions of
\begin{equation}\label{eq0 2.1}
 \arraycolsep=1pt\left\{
\begin{array}{lll}
 \displaystyle (-\Delta)^s u_i = f_i\quad
   &{\rm in}\quad  O\setminus \{0\},\\[2mm]
 \phantom{ (-\Delta)^s     }
 \displaystyle  u_i= h_i\quad  &{\rm   in}\ \ \R^N\setminus  O
 \end{array}\right.
\end{equation}
and
$$\limsup_{|x|\to0^+}u_1(x)|x|^{N-2s}\leq \liminf_{|x|\to0^+}u_2(x)|x|^{N-2s}.$$
  If $f_1\le f_2$ in $O\setminus \{0\}$ and $h_1\leq h_2$ in $\R^N\setminus  O$, then
$$u_1\le u_2\quad{\rm in}\quad O\setminus \{0\}.$$

\end{lemma}
{\bf Proof.}    Let $u=u_1-u_2$ and then
$$(-\Delta)^s u \le 0\quad {\rm in}\ \, O\setminus\{0\}\qquad {\rm and}\quad \limsup_{|x|\to0^+}u(x)\Phi_{s }^{-1}(x)\leq0,$$
 then for any $\epsilon>0$, there exists $r_\epsilon>0$ converging to zero as $\epsilon\to0$ such that
 $$u\le \epsilon \Phi_{s}\quad{\rm in}\quad \overline{B_{r_\epsilon} }\setminus\{0\},$$
 where $\Phi_{s}(x)=|x|^{2s-N}$ is the fundamental solution of $(-\Delta)^s$ in $\R^N$.

We see that
$$u=0<\epsilon \Phi_{s} \quad{\rm in}\ \ \R^N\setminus O,$$
then we have that
 $u\le \epsilon \Phi_s $ in $O\setminus\{0\}. $
By letting $\epsilon$ go to zero, we have that $u\le 0$  in $O\setminus\{0\}$.
  \hfill$\Box$\medskip

\noindent{\bf Proof of Theorem \ref{teo 2.1}. } By contradiction, we assume that $u$ is a positive solution of (\ref{eq 2.2-0}). Let $f_n(x)=(1- \eta_0(nx))f(x)$, where
$\eta_0:\R^N\to [0,1]$ is a radially symmetric,  smooth function such that
$$\eta_0=0\ \ {\rm in}\ \, \R^N\setminus B_{2}  \quad{\rm and}\quad \eta_0=1 \ \ {\rm in}\ \,  B_{1}.  $$
Since $f_n$ is H\"older continuous, it follows by \cite[Theorem 4.1]{CW} that
 \begin{equation}\label{eq 2.1-n}
\left\{\arraycolsep=1pt
\begin{array}{lll}
(-\Delta)^s  v=  f_n\qquad      {\rm in}\ \,  \Omega\setminus\{0\}, \\[2mm]
\qquad \quad  v=0\qquad   \,  \ {\rm in}\ \,   \R^N\setminus \Omega,\\[2mm]
\displaystyle \lim_{|x|\to0^+}v(x)|x|^{N-2s}=0
  \end{array}
 \right.
\end{equation}
  has unique solution $v_n$ and
  $$0<v_{n}\leq    u\quad {\rm in}\ \, \Omega\setminus\{0\} $$
   by Lemma \ref{cr hp} and the assumption that $h$ is non-negative. 
  
By the stability results \cite[Theorem 2.4]{CFQ}  and the regularity result  \cite[Theorem 2.1]{CFQ},       the limit $\{v_n\}_n$ exists, denoting $v_\infty$, is a positive classical solution of
  \begin{equation}\label{eq 2.1-l}
\left\{\arraycolsep=1pt
\begin{array}{lll}
(-\Delta)^s  v=  f\quad \     {\rm in}\ \  \Omega\setminus\{0\},\\[2mm]
\qquad \quad v=0\quad   \  {\rm in}\ \  \R^N\setminus \Omega,
  \end{array}
 \right.
\end{equation}
then we obtain a contradiction from  Theorem \ref{theorem-C} part $(ii)$.
  \hfill$\Box$ \smallskip

 \subsection{ Estimates }

For the calculations of  the fractional Laplacian for the involved logarithmic functions, we need the following estimates.  Recall that  for $\tau\in(-N,2s)$
\begin{equation}
  \label{eq:fractional-power}
 \cC_\s(\tau) = 2^{2\s} \frac{\Gamma(\frac{N+\tau}{2})\Gamma(\frac{2\s-\tau}{2})}{\Gamma(-\frac{\tau}{2})\Gamma(\frac{N-2\s+\tau}{2})}
\end{equation}
 and then
$$\cC_\s(\tau)=0$$
has two zero  points $0, 2s-N$.  Moreover, there holds
\begin{equation}
  \label{eq:fractional-power-1}
 (-\Delta)^\s |\cdot|^\tau = \cC_\s(\tau)|\cdot|^{\tau-2\s} \quad\  {\rm in}\ \  \mathcal{S}'(\R^N)
\end{equation}
and
$$
\cC_\s(\tau)  = -\frac{c_{N,\s}}2 \int_{\R^N}\frac{|e_1+z|^{\tau}+|e_1-z|^\tau-2}{|z|^{N+2\s}}\,dz,
$$
 where $e_1=(1,0,\cdots,0)\in\R^N$.

\begin{lemma}\cite[Lemma 2.3]{CW}
\label{c-function}
The function $\cC_s$, defined in (\ref{eq:fractional-power}),
is strictly concave and uniquely maximized at the point $\frac{2\s-N}{2}$ with the maximal value $2^{2\s}  \frac{\Gamma^2(\frac{N+2\s}4)}{\Gamma^2(\frac{N-2\s}{4})}.$

Moreover,
\begin{equation}
  \label{eq:c-symmetry}
\cC_\s(\tau)= \cC_\s(2\s-N-\tau) \qquad  {\rm for} \ \, \tau \in (-N,2\s)
\end{equation}
and
\begin{equation}
  \label{eq:c-asymptotics}
\lim_{\tau \to -N}\cC_\s(\tau) = \lim_{\tau \to 2 \s}\cC_\s(\tau)=-\infty.
\end{equation}
\end{lemma}

Note that
 $$ \cC_\s(\tau)>0\quad{\rm  for}\ \, \tau\in (2s-N, 0), $$
 $$ \cC_\s(\tau) <0\quad{\rm  for}\ \, \tau\in (0,2s)\cup (-N,2s-N)$$
and
 $$ \cC_\s'(2s-N)=-\cC_\s'(0)>0\quad {\rm and}\quad \cC_\s''(2s-N)=\cC_\s''(0)<0. $$

For $m\in\R$  denote  $ v_{m}$ be a smooth,  radially symmetric   function   such that
 \begin{equation}\label{cl-6-1}
 v_{ m}(x)= (-\ln|x|)^m\quad\ {\rm for}\ \ 0<|x|<\frac1{e^2}\quad {\rm and}\quad v_{ m}(x) =0\quad {\rm for }\ |x|>1.
 \end{equation}
 Moreover,  $v_m$ is non-increasing  in $|x|$ if $m>0$.
Denote
 \begin{equation}\label{cl-6}
 w_{m}(x)=|x|^{2s-N} v_m(x) \ \ {\rm for}\ \,  x\in \R^N\setminus\{0\}.
 \end{equation}

 \begin{proposition}\label{pr 2.1-a}
Let $m\not=0$ and
 \begin{equation}\label{cl-6-2}
 \cB_{m}=\cC_s'(0)m,\quad \cD_{m}=\cC_s''(0)\frac{m(m-1)}2,
 \end{equation}
then $\cB_{m}>0$ for $m<0$, and there exist $r_0\in(0,\frac1{e^2}]$ and $c_5>0$ such that
 \begin{eqnarray}
  \Big| (-\Delta)^\s w_{m}(x)  -\cB_{m}  |x|^{-N}(-\ln |x|)^{m-1}-\cD_{m}|x|^{-N}(-\ln|x|)^{m-2}\Big| \nonumber
 \leq   c_5 |x|^{-N}(-\ln |x|)^{m-3}.\label{cl-4}
\end{eqnarray}
\end{proposition}

In order to get precise estimates of  $(-\Delta)^\s w_{m}$, we use the  fact  that  for $x\in B_{\frac1{e^2}}\setminus\{0\}$,
 \begin{eqnarray*}
 (-\Delta)^\s w_{m}(x)   &= &v_m(x) (-\Delta)^\s |x|^{2s-N} + |x|^{2s-N}(-\Delta)^\s  v_m(x) 
\\[2mm]&& +
c_{s,N}  \int_{\R^N}\frac{(|x|^{2s-N}-|y|^{2s-N})(v_m(x)-v_m(y))}{|x-y|^{N+2s}}dy
\\[2mm]&=&F_m(x)+E_m(x),
\end{eqnarray*}
where $(-\Delta)^\s |x|^{2s-N}=0$ in $\R^N\setminus\{0\}$,
$$ E_m(x)= c_{s,N}  \int_{\R^N}\frac{(|x|^{2s-N}-|y|^{2s-N})(v_m(x)-v_m(y))}{|x-y|^{N+2s}}dy$$
and
$$F_m(x)=  |x|^{2s-N}(-\Delta)^\s  v_m(x).$$
In order to get   estimate (\ref{cl-4}), we need the following lemmas.

\begin{lemma}\label{lm 2.1-c}
Let $m<0$, then there exist $r_0\in(0,\frac1{e^2}]$ and $c_6>0$ such that  for $0<|x|<r_0$
 \begin{equation}\label{cl-2}
 \Big|E_m(x)-2\cB_{m}  |x|^{-N}(-\ln |x|)^{m-1}\Big|\leq  c_6 |x|^{s-N}(-\ln |x|)^{m}.
  \end{equation}
\end{lemma}
{\bf Proof.}
A direct computation shows that
 \begin{eqnarray*}
\cC_\s'(\tau)  &=& -\frac{c_{N,\s}}2 \int_{\R^N}\frac{|e_1+z|^{\tau}\ln|e_1+z|+|e_1-z|^\tau\ln|e_1-z| }{|z|^{N+2\s}}\,dz
\\[1.5mm]&=&-\frac{c_{N,\s}}2\Big( \int_{\R^N}\frac{\big(|e_1+z|^{\tau}-1\big)\ln|e_1+z|  }{|z|^{N+2\s}}\,dz
\\[1.5mm]&&   +\int_{\R^N}\frac{\big(|e_1-z|^{\tau}-1\big)\ln|e_1-z|  }{|z|^{N+2\s}}\,dz
+\int_{\R^N}\frac{ \ln|e_1-z| +\ln|e_1+z|  }{|z|^{N+2\s}}\,dz\Big)
\\[1.5mm]&=&c_{N,\s}\int_{\R^N}\frac{\big(1-|z|^{\tau}\big)   \big(-\ln |z| \big)   }{|e_1-z|^{N+2s}}dz
-\frac{c_{N,\s}}2 \int_{\R^N}\frac{ \ln|e_1-z| +\ln|e_1+z|  }{|z|^{N+2\s}}\,dz
\end{eqnarray*}
and
 \begin{eqnarray*}
\cC_\s''(\tau)  &=& -\frac{c_{N,\s}}2 \int_{\R^N}\frac{|e_1+z|^{\tau}(\ln|e_1+z|)^2+|e_1-z|^\tau(\ln|e_1-z|)^2 }{|z|^{N+2\s}}\,dz
\\&=&c_{N,\s}\int_{\R^N}\frac{\big(1-|z|^{\tau}\big)   \big(-\ln |z| \big)^2   }{|e_1-z|^{N+2s}}dz-
c_{N,\s}\int_{\R^N}\frac{ \big(-\ln |z| \big)^2   }{|e_1-z|^{N+2s}}dz.
\end{eqnarray*}
Henceforth
$$\cC_\s'(0)=-\frac{c_{N,\s}}2 \int_{\R^N}\frac{ \ln|e_1-z| +\ln|e_1+z|  }{|z|^{N+2\s}}\,dz$$   and
$$ \cC_\s''(0)=-c_{N,\s}\int_{\R^N}\frac{ \big(-\ln |z| \big)^2   }{|e_1-z|^{N+2s}}dz. $$
Therefore, we obtain some important equalities:
 \begin{eqnarray*}
&&\cC_\s'(\tau) -\cC_\s'(0)=c_{N,\s}\int_{\R^N}\frac{ (1-|z|^{\tau} )    (-\ln |z|  )   }{|e_1-z|^{N+2s}}dz, \\[2mm]
\\&& \cC_\s''(\tau) -\cC_\s''(0)=c_{N,\s}\int_{\R^N}\frac{ (1-|z|^{\tau} )    (-\ln |z|  )^2   }{|e_1-z|^{N+2s}}dz
\end{eqnarray*}
and
$$c_{N,\s}\int_{\R^N}\frac{\big(1-|z|^{2s-N}\big)   \big(-\ln |z| \big)   }{|e_1-z|^{N+2s}}dz=\cC_\s'(2s-N) -\cC_\s'(0)=-2\cC_\s'(0)>0,$$
 \begin{equation}\label{cl-2-ss}
 c_{N,\s}\int_{\R^N}\frac{ (1-|z|^{2s-N} )    (-\ln |z|  )^2   }{|e_1-z|^{N+2s}}dz=\cC_\s''(2s-N) -\cC_\s''(0)=0.
 \end{equation}

Let $0<|x|<r_0$ with $r_0>0$ small enough, we see that
\begin{eqnarray*}
&&\int_{B_{\sqrt{|x|}}(x)}\frac{(|x|^{2s-N}-|y|^{2s-N})(v_m(x)-v_m(y))}{|x-y|^{N+2s}}dy
 \\&=& |x|^{-N}(-\ln |x|)^m \int_{B_\frac1{\sqrt{|x|}}(e_1)}\frac{\big(1-|z|^{2s-N}\big)\Big(1-\big(1+\frac{-\ln |z| }{-\ln|x|}\big)^m \Big)}{|e_x-z|^{N+2s}}dz,
 \end{eqnarray*}
 where
  $\big|\frac{-\ln |z|}{-\ln|x|}\big|\leq \frac12$ for $|z|<\frac1{\sqrt{|x|}}$ and
 $$ 1-\big(1+\frac{-\ln |z|}{-\ln|x|}\big)^m=-m\frac{-\ln |z|}{-\ln|x|}+O\Big(\big(\frac{-\ln |z|}{-\ln|x|}\big)^2\Big).  $$
Thus, we see that
 \begin{eqnarray*}
 && \int_{B_{\sqrt{|x|}}(x)}\frac{(|x|^{2s-N}-|y|^{2s-N})(v_m(x)-v_m(y))}{|x-y|^{N+2s}}dy
  \\[2mm]&=&-m |x|^{-N}(-\ln |x|)^{m-1} \int_{B_\frac1{\sqrt{|x|}}(e_1)}\frac{\big(1-|z|^{2s-N}\big)   \big(-\ln |z| \big)   }{|e_1-z|^{N+2s}}dz
  \\&&\qquad     +  |x|^{-N}(-\ln |x|)^{m-2} O \Big(\int_{B_\frac1{\sqrt{|x|}}(e_1)}\frac{\big(1-|z|^{2s-N}\big)   \big(-\ln |z| \big)^2   }{|e_1-z|^{N+2s}}dz\Big),
 \end{eqnarray*}
 where    \begin{eqnarray*}
&& \Big|\frac{-2\cC_\s'(0)}{c_{N,s}}-\int_{B_\frac1{\sqrt{|x|}}(e_1)}\frac{\big(1-|z|^{2s-N}\big)   \big(-\ln |z| \big)    }{|e_1-z|^{N+2s}}dz \Big|
\\[1.5mm] &=& \Big|  \int_{\R^N\setminus B_\frac1{\sqrt{|x|}}(e_1)}\frac{\big(1-|z|^{2s-N}\big)   \big(-\ln |z| \big)    }{|e_1-z|^{N+2s}}dz\Big|
 \\[1.5mm] &\leq &\int_{\R^N\setminus B_\frac1{\sqrt{|x|}}(e_1)}\frac{    -\ln |z|     }{|e_1-z|^{N+2s}}dz\\[1.5mm] &\leq &  c_7 \big(- \ln |x| \big)  |x|^s
  \end{eqnarray*}
  and
 \begin{eqnarray*}
 &&\Big|\int_{B_\frac1{\sqrt{|x|}}(e_1)}\frac{\big(1-|z|^{2s-N}\big)   \big(-\ln |z| \big)^2   }{|e_1-z|^{N+2s}}dz\Big|\\ [1.5mm] &=& \Big| \int_{\R^N}\frac{\big(1-|z|^{2s-N}\big)   \big(-\ln |z| \big)^2   }{|e_1-z|^{N+2s}}dz-\int_{\R^N\setminus B_\frac1{\sqrt{|x|}}(e_1)}\frac{\big(1-|z|^{2s-N}\big)   \big(-\ln |z| \big)^2   }{|e_1-z|^{N+2s}}dz\Big|
 \\[1.5mm] &\leq &\int_{\R^N\setminus B_\frac1{\sqrt{|x|}}(e_1)}\frac{   \big(-\ln |z| \big)^2   }{|e_1-z|^{N+2s}}dz\\[1.5mm]&\leq &  c_7 \big(\frac12 \ln |x| \big)^2 |x|^s
  \end{eqnarray*}
 by using (\ref{cl-2-ss}).

 As a consequence, we obtain that
 \begin{eqnarray}
&& \nonumber \int_{B_{\sqrt{|x|}}(x)}\frac{(|x|^{2s-N}-|y|^{2s-N})(v_m(x)-v_m(y))}{|x-y|^{N+2s}}dy
\\[1.5mm] &=& 2\cB_m  |x|^{-N}(-\ln |x|)^{m-1} \Big(1  +O\big(|x|^s(-\ln |x|)\big) \Big).\label{cl-1}
\end{eqnarray}

On the other hand, we see that
  \begin{eqnarray*}
&& \Big| \int_{\R^N\setminus B_{\sqrt{|x|}}(x)}\frac{(|x|^{2s-N}-|y|^{2s-N})(v_m(x)-v_m(y))}{|x-y|^{N+2s}}dy\Big|
\\[1.5mm] &=& |x|^{2s-N} \big((-\ln |x|)^{m}+1\big) \int_{\R^N\setminus B_{\sqrt{|x|}}(x)}\frac{ 1 }{|x-y|^{N+2s}}dy
\\[1.5mm] &\leq & c_8|x|^{s-N} \big((-\ln |x|)^{m}+1\big),
 \end{eqnarray*}
where  $c_8>0$ is independent of $|x|$ and
$$v_m(y)\leq v_m(x)\quad{\rm for}\ \, |y|>|x|.$$
Together with (\ref{cl-1}), we obtain (\ref{cl-2}).
Then we complete the proof.\hfill$\Box$\medskip

\begin{lemma}\label{lm 2.1-d}
There exist $r_0\in(0,\frac1{e^2}]$ and $c_9>0$ such that    for $0<|x|<r_0$
 \begin{equation}\label{cl-3}
 \begin{array}{lll}
 \Big|F_m(x)+\cB_{m}  |x|^{-N}(-\ln |x|)^{m-1}-\cD_m|x|^{-N}(-\ln |x|)^{m-2}\Big|
\leq  c_9 |x|^{-N}(-\ln |x|)^{m-3}.
  \end{array}
  \end{equation}
\end{lemma}
\noindent {\bf Proof. } For $x\in B_{\frac1{2e}}\setminus\{0\}$, we see that
 \begin{eqnarray*}
(-\Delta)^\s v_{m}(x)  &= &\frac{c_{N,s}}{2}\int_{B_{\sqrt{|x|}}} \frac{ 2(-\ln|x|)^m-(-\ln|x+y|)^m-(-\ln|x-y|)^m  }{|y|^{N+2s}}dy\\[1.5mm]&&+c_{N,s}\int_{\R^N\setminus B_{\sqrt{|x|}}}\frac{  (-\ln|x|)^m-v_m(y)  }{|y|^{N+2s}}dy,
\end{eqnarray*}
where
 \begin{eqnarray*}
 0&<&c_{N,s}\int_{\R^N\setminus B_{\sqrt{|x|}}}\frac{  (-\ln|x|)^m-v_m(x+y)  }{|y|^{N+2s}}dy
 \\[1.5mm]&\leq& c_{N,s}\int_{\R^N\setminus B_{\sqrt{|x|}}}\frac{  (-\ln|x|)^m   }{|y|^{N+2s}}dy
   \\[1.5mm]&\leq&   c_{10}|x|^{-s} (-\ln|x|)^m.
\end{eqnarray*}
 Moreover, there holds
 \begin{eqnarray*}
 &&\frac{c_{N,s}}{2}\int_{B_{\sqrt{|x|}}} \frac{ 2(-\ln|x|)^m-(-\ln|x+y|)^m-(-\ln|x-y|)^m  }{|y|^{N+2s}}dy
 \\[1.5mm] &=& \frac{c_{N,s}}{2} |x|^{-2s}(-\ln|x|)^m  \int_{ B_{\frac1{\sqrt{|x|}}}}\frac{2-(1+\frac{-\ln|z+e_1|}{-\ln|x|})^m  -(1+\frac{-\ln|z-e_1|}{-\ln|x|})^m  }{|z|^{N+2s}}dz
 \\[1.5mm]&=&  \frac{c_{N,s} m}{2}   |x|^{-2s}(-\ln|x|)^{m-1} \int_{ B_{\frac1{\sqrt{|x|}}}}\frac{ \ln|e_1-z| +\ln|e_1+z|  }{|z|^{N+2\s}}\,dz
 \\[1.5mm]&&\quad  -  \frac{c_{N,s} m(m-1)}{4}   |x|^{-2s}(-\ln|x|)^{m-2}
 \int_{ B_{\frac1{\sqrt{|x|}}}}\frac{ (\ln|e_1-z|)^2 +(\ln|e_1+z|)^2  }{|z|^{N+2\s}}\,dz
  \\[1.5mm]&&\quad  +|x|^{-2s}   (-\ln|x|)^{m-3} O\Big(\int_{\R^N}\frac{ (\ln|e_1-z|)^3+ (\ln|e_1-z|)^3   }{|z|^{N+2\s}}\,dz\Big)
 \\[1.5mm]&=&m |x|^{-2s}(-\ln|x|)^{m-1}\Big(-\cC_s'(0)- c_{N,s} \int_{ \R^N\setminus B_{\frac1{\sqrt{|x|}}}}\frac{ \ln|e_1-z|    }{|z|^{N+2\s}}\,dz \Big)
 \\[1.5mm]&&\quad +\frac{m(m-1)}2 |x|^{-2s}(-\ln|x|)^{m-2}\Big(\cC_s''(0)- c_{N,s} \int_{ \R^N\setminus B_{\frac1{\sqrt{|x|}}}}\frac{ (\ln|e_1-z|)^2    }{|z|^{N+2\s}}\,dz \Big)
 \\[1.5mm]&&\quad  +|x|^{-2s}   (-\ln|x|)^{m-3} O(1)
 \\[1.5mm]&=&|x|^{-2s} (-\ln|x|)^{m-1}\Big( -m\cC_s'(0)+\frac{m(m-1)}2\cC_s''(0)(-\ln|x|)^{-1} + (-\ln|x|)^{-2} O(1)\Big),
\end{eqnarray*}
where we used  $|\frac{-\ln|z\pm e_1|}{-\ln|x|}|\ll 1$ when $|x|$ small enough and
$$(1+t)^m=1+mt+\frac{m(m-1)}{2}t^2+O(t^3).$$

As a consequence, we conclue that
 \begin{eqnarray*}
&& \Big|(-\Delta)^\s v_{m}(x) +m\cC_s'(0)|x|^{-2s}(-\ln|x|)^{m-1}
-\frac{m(m-1)}2\cC_s''(0)|x|^{-2s}(-\ln|x|)^{m-2}\Big|
\\[1.5mm]&\leq& c_{11}|x|^{-2s}(-\ln|x|)^{m-3},
\end{eqnarray*}
which implies (\ref{cl-3}).
We complete the proof.\hfill$\Box$\medskip

 \noindent{\bf Proof of Proposition \ref{pr 2.1-a}. } It follows by Lemma \ref{lm 2.1-c} and
Lemma \ref{lm 2.1-d} directly. \hfill$\Box$\medskip

From the proof of Lemma \ref{lm 2.1-d}, we have the following corollary.
 \begin{corollary}\label{cr 2.1-a}
Let $m\not=0$
then   there exist $r_0\in(0,\frac1{e^2})$ and $c_{12}>0$ such that for $0<|x|<r_0$
 \begin{equation}\label{cl-4-l}
 \begin{array}{lll}
\displaystyle  \Big| (-\Delta)^\s v_{m}(x)  +\cB_{m}  |x|^{-2s}(-\ln |x|)^{m-1}-\cD_{m}|x|^{-2s}(-\ln|x|)^{m-2}\Big|
 \leq   c_{12}|x|^{-2s}(-\ln |x|)^{m-3}.
  \end{array}
  \end{equation}
\end{corollary}

More generally, for $\tau\in(-N,2s),\, m\in\R$,  let
$$w_{\tau, m}(x)=|x|^{\tau} v_m(x) \ \ {\rm for}\ \,  x\in \R^N\setminus\{0\}.$$
The same calculation implies that
\begin{corollary}\label{cr 2.1-a}
Let $\tau\in(-N,2s),\, m\in\R$, $|\tau|+|m|\not=0$ and
$$\cB_{\tau, m}= \cC_s'(\tau) m,$$
then  there exist $r_0\in(0,\frac1{e^2})$ and $c_{13}>0$ such that for $0<|x|<r_0$
 \begin{equation}\label{cl-5}
  \begin{array}{lll}
\displaystyle \Big| (-\Delta)^\s w_{\tau, m}(x) -\cC_s(\tau) w_{\tau, m}(x) |x|^{-2s} -\cB_{\tau,m}  |x|^{\tau-2s}(-\ln |x|)^{m-1}\Big|
 \leq   c_{13}|x|^{\tau-2s}(-\ln |x|)^{m-2}.
  \end{array}
   \end{equation}
\end{corollary}

\subsection{Blow-up rate  estimates}
The following estimates play an important role in our classification of isolated singularities of (\ref{eq 1.1}).   In what follows, we always let $r_0\in(0,\frac1{e^2})$ be from Proposition \ref{pr 2.1-a}.  

   \begin{lemma}\label{lm 2.2-lower} Let $g\in C^\beta_{loc}(\Omega\setminus\{0\})$,  with $\beta\in(0,1)$,  be a   nonnegative function such that  for some  $m<0 $  
 $$ g(x)\geq |x|^{-N} (-\ln|x|)^{m-1}\quad{\rm in}\ \,  B_{r_0}\setminus\{0\}.$$
 Let  $u_g$ be a positive solution of problem
 \begin{equation}\label{eq 2.1-hom}
 \left\{\arraycolsep=1pt
\begin{array}{lll}
 (-\Delta)^\s  u=  g\quad\ \
   {\rm in}\ \,  \Omega\setminus \{0\},\\[2mm]
\qquad \ \ u\geq  0\quad \ \  {\rm   in}\ \,  \R^N\setminus \Omega.
  \end{array}
 \right.
 \end{equation}

Then   there exists $r\in(0,r_0]$    such that
$$u_g(x)\geq \cB_m   |x|^{2s-N}(-\ln|x|)^{m}\quad{\rm in}\ \  B_r\setminus\{0\}.$$

\end{lemma}
\noindent{\bf Proof. }  Recall that $w_m$ is defined in (\ref{cl-6}), has compact support in $\bar B_1$ and
$$w_m(x)=|x|^{2s-N}(-\ln|x|)^{m}\quad  {\rm in}\ \ B_{\frac1{e^2}}\setminus\{0\}$$
by  Proposition \ref{pr 2.1-a}, there exists $r_1\in(0,r_0]$ such that
 \begin{eqnarray*}
(-\Delta)^\s w_m &\leq&  \cB_m |x|^{-N}(-\ln|x|)^{m-1}+\cD_m |x|^{-N}(-\ln|x|)^{m-2}
\\[1.5mm]&&\quad +O(1)|x|^{-N}(-\ln|x|)^{m-3}
\\[1.5mm]&\leq &   \cB_m  |x|^{-N}(-\ln|x|)^{m-1} \quad{\rm for}\ x\in B_{r_1}\setminus\{0\},
\end{eqnarray*}
where   $\cB_m>0$ and $\cD_m<0$ for $m<0$.

Note that $u_g$ is positive and continuous in $\Omega\setminus\{0\}$,
then  $u_g\geq  0$ in $\R^N\setminus B_{r_0}$.   By the lower bound of $g$ there exists $t_0=  \cB_m >0$ such that
 \begin{eqnarray*}
 (-\Delta)^\s t_0u_g\geq t_0g(x)  &\geq& t_0 |x|^{-N}(-\ln|x|)^{m-1}
 \\[1.5mm]&\geq &(-\Delta)^\s w_m
 \\[1.5mm]&= & (-\Delta)^\s \big(w_m-w_m(r_1)\big)\quad{\rm in}\quad B_{r_1}\setminus\{0\}
 \end{eqnarray*}
 and
 $$\liminf_{|x|\to0^+}u_g\Phi_s^{-1}(x)\geq0= \lim_{|x|\to0^+}\big(w_m-w_m(r_1)\big) \Phi_s^{-1}(x),$$
 $$ t_0u_g\geq 0\geq \big(w_m-w_m(r_1)\big)\quad{\rm in}\ \ \R^N\setminus B_{r_1}.$$
Then  by Lemma \ref{cr hp}, we have that
$$ t_0u_g\geq  w_m-w_m(r_1) \quad{\rm in}\quad \Omega\setminus\{0\}.$$
This completes the proof.  \hfill$\Box$\medskip

\begin{lemma}\label{lm 2.2-upper}  Let $g\in C^\beta_{loc}(\Omega\setminus\{0\})$,  with $\beta\in(0,1)$,  be a   nonnegative function such that    $$ g(x)\leq |x|^{-N} (-\ln|x|)^{m-1}\quad{\rm in}\ \,  B_r\setminus\{0\}.$$
 Let $h\in C(B_{2R_0})\cap L^1_s(\R^N)$ be a nonnegative function such that 
 $$h=0\quad {\rm in} \ B_{ \frac12  }.$$

 Let  $u_g$ be a positive solution of the problem:
 \begin{equation}\label{eq 2.1-hom-}
 \left\{\arraycolsep=1pt
\begin{array}{lll}
(-\Delta)^\s  u= g\quad\ \
   {\rm in}\ \,  \Omega\setminus \{0\},\\[2mm]
\qquad \ \ u= h\quad \ \  {\rm   in}\ \,  \R^N\setminus \Omega,\\[2mm]
 \displaystyle\lim_{|x|\to0^+} u(x)|x|^{N-2s} =0.
  \end{array}
 \right.
 \end{equation}
Then for any $\epsilon\in(0,\cB_m)$,  there exists $r\in(0,r_0]$   such that
$$u_g(x)\leq ( \cB_m+\epsilon)  |x|^{2s-N}(-\ln|x|)^{m}\quad{\rm in}\ \  B_r\setminus\{0\}.$$

\end{lemma}
\noindent{\bf Proof. }  Let $U_0(x)=-(-\Delta)^s h(x)\quad {\rm for }\ x\in B_{r_0}$. A direct computation shows that for $x\in B_{r_0}$
\begin{eqnarray}
 0<U_0(x) &=&  c_{N,s} \int_{\R^N\setminus B_{\frac12}}\frac{h(y)}{|x-y|^{N+2s}}dy\nonumber
 \\ &\leq &    c_{N,s}   \int_{\R^N\setminus B_{\frac12}}\frac{h(y)}{(|y|-r_0)^{N+2s}}dy
 \label{est 2.1}  \end{eqnarray}
by using the fact that
$$|x-y|\geq |y|-r_0\quad {\rm for}\ \, |x|\leq r\ \, {\rm and}\ \, |y|\geq \frac12,$$
and taking into account that
$$w_m(x)=|x|^{2s-N}(-\ln|x|)^{m}\quad  {\rm in}\ \ B_{\frac1{e^2}}\setminus\{0\}$$
and  that by Proposition \ref{pr 2.1-a}, for given $\epsilon\in(0, \cB_m)$, there exists $r_2\in(0,r_0]$ such that
 \begin{eqnarray*}
(-\Delta)^\s w_m &\geq&  \cB_m |x|^{-N}(-\ln|x|)^{m-1}-c_{13} |x|^{-N}(-\ln|x|)^{m-2}
\\[1.5mm]&\geq & ( \cB_m-\epsilon) |x|^{-N}(-\ln|x|)^{m-1} \quad{\rm for}\ x\in B_{r_2}\setminus\{0\}.
\end{eqnarray*}

Note that $u_g$ is positive and continuous in $\Omega\setminus\{0\}$,
then  $u_g\leq \varrho_0$ in $\Omega\setminus B_{r_2}$ for some $\varrho_0>0$.   By the upper bound of $g$ there exists $t_1=( \cB_m-2\epsilon)>0$ such that
 \begin{eqnarray*}
 (-\Delta)^\s t_1(u_g-\varrho_0-h)&\leq& t_1(g(x)-U_0)
 \\[1.5mm]&\leq& ( \cB_m-\epsilon) |x|^{-N}(-\ln|x|)^{m-1}
 \\[1.5mm]&\leq &(-\Delta)^\s w_m
 \\[1.5mm]&= &(-\Delta)^\s \big(w_m-w_m(r_1)\big)\quad{\rm in}\quad B_{r_2}\setminus\{0\},
 \end{eqnarray*}
 where we used that $U_0$ is bounded in $B_{r_2}$.
Moreover, we have that
 $$\liminf_{|x|\to0^+}\big((u_g(x)-\varrho_0-h(x))\Phi_s^{-1}(x)\big)=0= \lim_{|x|\to0^+} w_m \Phi_s^{-1}(x) $$
 and
 $$  t_1\big(u_g-\varrho_0-h\big)=t_1\big(u_g-\varrho_0\big)\leq 0\leq  w_m \quad{\rm in}\ \ \R^N\setminus B_{r_2}.$$
Then  by Lemma \ref{cr hp}, we have that
$$  t_1\big(u_g-\varrho_0-h\big) \leq  w_m  \quad{\rm in}\quad \Omega\setminus\{0\}.$$
This completes the proof.   \hfill$\Box$\medskip

\begin{lemma}\label{teo 2.3}  Let $g\in C^\beta_{loc}(\Omega\setminus\{0\})$,  with $\beta\in(0,1)$,  be a   nonnegative function such that  there exists $\tau\in(2s-N,+\infty) $  such that
 $$ g(x)\leq  |x|^{\tau-2s}\quad{\rm in}\ \, B_{r_0}\setminus\{0\}.$$
 Let $h\in C(B_{2R_0})\cap L^1_s(\R^N)$ be a nonnegative function such that 
 $$h=0\quad {\rm in} \ B_{\frac12}.$$

 Let  $u_g$ be a positive solution of problem
 \begin{equation}\label{eq 2.1-hom-}
 \left\{\arraycolsep=1pt
\begin{array}{lll}
(-\Delta)^\s  u\leq g\quad\ \
   {\rm in}\ \,  \Omega\setminus \{0\},\\[2mm]
\qquad \ \ u= h\quad \ \  {\rm   in}\ \,  \R^N\setminus \Omega,\\[2mm]
 \displaystyle\lim_{|x|\to0^+} u(x)|x|^{N-2s} =0.
  \end{array}
 \right.
 \end{equation}
Then 

$(i)$ for $\tau\in(2s-N ,0)$ there exists $c_{14}>0$  such that
$$u_g(x)\leq c_{14}|x|^\tau\quad{\rm in}\ \  \Omega\setminus\{0\};$$

$(ii)$ for $\tau\in(0,2s)$, there exists $c_{15}>0$  such that
$$u_g(x)\leq c_{15} \quad{\rm in}\ \  \Omega\setminus\{0\}.$$

\end{lemma}
\noindent{\bf Proof.}   Recall that $U_0(x)=-(-\Delta)^s h(x)\quad {\rm for }\ x\in B_{r_0}$, which is bounded by  (\ref{est 2.1}).

$(i)$ For $\tau\in(2s-N ,0)$, we have that
$$(-\Delta)^s |x|^{\tau}=\cC_s(\tau) |x|^{\tau-2s}\quad  {\rm in}\ \  \R^N\setminus\{0\},$$
where $\cC_s(\tau)>0$.

Note that  there exists $t_3>0$ such that
 $$(-\Delta)^s t_3(u_g-h)\leq t_3(g(x)+U_0) \leq  \cC_s(\tau) |x|^{\tau-2s} =(-\Delta)^s |x|^{\tau} \quad{\rm in}\quad \Omega\setminus\{0\}  $$
 and
 $$\lim_{|x|\to0}u_g(x)|x|^{N-2s}=0,\qquad  u_g= h \quad{\rm in}\ \ \R^N\setminus \Omega.$$
Then  by Lemma \ref{cr hp}, we have that
$$ t_3u_g(x)\leq  |x|^{\tau} \quad{\rm for}\ \, x\in \Omega\setminus\{0\}.$$

$(ii)$   Take
$$\tilde w_3(x)=(2R_0)^{\tau }  -|x|^{\tau}\quad  {\rm in}\ \  \R^N\setminus\{0\}$$
and
$$w_3(x)=(2R_0)^{\tau } -|x|^{\tau}\quad  {\rm in}\ \  B_{2R_0}\setminus\{0\},
\qquad w_3(x)=0\quad{\rm in}\ \ \R^N \setminus B_{R_0}.$$
Direct computation shows that for $x\in B_{R_0}\setminus\{0\}$
\begin{eqnarray*}
(-\Delta)^s w_3(x)  &=&(-\Delta)^s \tilde  w_3(x) + (-\Delta)^s (w_3-\tilde  w_3)(x)
\\[1.5mm]&\geq &  -\cC_{s}(\tau) |x|^{\tau-2s}-c_{N,s}(2R_0)^{\tau }  \int_{\R^N\setminus B_{2R_0}} \frac{1 }{|x-y|^{N+2s}}dy
\\[1.5mm]&\geq &-\cC_{s}(\tau) |x|^{\tau-2s}-c_{N,s}(2R_0)^{\tau-2s} \int_{\R^N\setminus B_1} \frac{ 1 }{|\tilde e_x-z|^{N+2s}}dz
\\[1.5mm] &\geq &-\cC_{s}(\tau) |x|^{\tau-2s}-2^{N+\tau}c_{N,s}R_0^{\tau-2s}\int_{\R^N\setminus B_1} \frac{1}{|z|^{N+2s}}dz,
\end{eqnarray*}
where $-\cC_{s}(\tau)>0$ for $\tau\in(0,2s)$,  $\tilde e_x=\frac{x}{2R_0}$ and $|\tilde e_x-z|\leq \frac{|z|}2$.
 Thus, there exists $r\in (0,1)$ such that for $x\in B_r\setminus\{0\}$
 \begin{eqnarray*}
(-\Delta)^s w_3(x)   \geq  -\frac12\cC_{s}(\tau) |x|^{\tau-2s}.
\end{eqnarray*}

Consequently  there exists $t_4>0$ such that
 $$(-\Delta)^s t_4(u_g-h)\leq t_4(g +U_0) \leq  -\frac12\cC_{s}(\tau) |x|^{\tau-2s}\leq (-\Delta)^s w_3 \quad{\rm in}\quad B_r\setminus\{0\}  $$
 and
 $$w_3\geq (2^\tau-1) R_0^{\tau}\geq t_4(u_g-h)  \quad {\rm in}\ \, \overline{\Omega }\setminus B_r, $$
 since $\Omega\subset B_{R_0}$. 
 Together with
 $$\lim_{|x|\to0^+}u_g(x)|x|^{N-2s}=0,\qquad  u_g= h \quad{\rm in}\ \ \R^N\setminus B_r,$$
it implies by  Lemma \ref{cr hp}  that
$$ t_4(u_g-h) \leq  w_3 \quad{\rm for\ any}\ \, x\in  B_r\setminus\{0\}.$$ 

Therefore,  we obtain that
$$u_g(x)\leq c_{15} \quad{\rm for\ any}\ \, x\in  \Omega\setminus\{0\}. $$
This completes the proof. \hfill$\Box$

\setcounter{equation}{0}
\section{ Isolated singularity   }

 In this section, we provide  rough bounds for the isolated singular solution of (\ref{eq 1.1}).

 \begin{theorem}\label{teo 3-1}

  Let $u$ be a  positive solution of (\ref{eq 1.1}) verifying
  $$\lim_{|x|\to0^+} u(x)=+\infty,$$
  then
$${\rm \it  upper\ bound:}\qquad\qquad \limsup_{|x|\to0^+}u(x) |x|^{N-2s}<+\infty \qquad\qquad\qquad\qquad\qquad\qquad\qquad\qquad\qquad \qquad\qquad $$
and for any $\tau\in(2s-N,0)$
$${\rm \it  lower\ bound:}\qquad\qquad \lim_{|x|\to0^+}u(x) |x|^{-\tau}=+\infty. \qquad\qquad\qquad\qquad\qquad\qquad\qquad\qquad\qquad\qquad\qquad \quad $$

 \end{theorem}

 \subsection{ Upper bound}

 \begin{proposition}\label{pr 2.1}
 Let  $h\in C^\theta(B_{2R_0})\cap L^1_s(\R^N)$ with $\theta>2s$, 
 and 
  $u$ be a nonnegative solution of (\ref{eq 1.1}), then
  there exist $r_1\in(0,\frac1{e^2}]$ and $c_{16}>0$ such that
 \begin{equation}\label{ll 3.0}
 u(x)\leq c_{16} |x|^{2s-N} ,\quad\forall\, x\in B_{r_1}\setminus\{0\}.
 \end{equation}

 \end{proposition}

 In order to prove Proposition \ref{pr 2.1}, we need following lemma.

  \begin{lemma}\label{lm 3.1-1}
 Assume that  $h\in C^\theta(B_{2R_0})\cap L^1_s(\R^N)$ with $\theta>2s$, 
 and  $u$ is a nonnegative classical solution of (\ref{eq 1.1}) replaced $p^*$ by $p>1$.
  Then
$
 u^p\in L^1(\Omega,\rho dx)
$
and there exists a uniform $c_{17}>0$ independent of $u$ such that
\begin{equation}\label{lp}
\int_{\Omega} u^p \rho(x)^s dx<c_{17},
\end{equation}
where $\rho(x)={\rm dist}(x,\partial\Omega)$.
 \end{lemma}
\noindent{\bf Proof. }    Recall that $U_0(x)=-(-\Delta)^s  h(x)\quad {\rm for }\ x\in \Omega$, which is uniformly  bounded in $\Omega$.
 Let  $w=u-h$ in $\R^N$, then we have that
 $$(-\Delta)^s w=u^p+U_0\quad {\rm in}\ \, \Omega\setminus\{0\},\qquad w=0\quad {\rm in}\ \R^N\setminus\Omega$$
 From Theorem \ref{theorem-C}, we have that $u^p-U_0\in L^1(\Omega)$, so is $u$, thanks to the boundedness of $U_0$. Moreover, we have that
 $$\int_{\Omega}w   (-\Delta)^\s \xi \,dx = \int_{\Omega}(u^p+U_0) \xi\,   dx \qquad  {\rm for\ all}\ \,\xi \in  C^2_0(\Omega).   $$
 Let $(\lambda_1,\xi_1)$ be the the first eigenvalue and related positive eigenfunction of $(-\Delta)^s$  in $\Omega$ with zero Dirichlet boundary condition, i.e.
 $$(-\Delta)^\s \xi_1=\lambda_1\xi_1\quad  {\rm in}\ \, \Omega,\qquad \xi_1=0\quad {\rm in}\ \R^N\setminus\Omega.$$
  The existence and properties could see \cite[Proposition 9]{SV}. In fact, we have that  $\lambda_1>0$, $\xi_1$ is positive,  $\xi_1\in C^{2s}_{loc}(\Omega)\cap C^s(\R^N)$
and $\xi(x)\sim \rho^s(x)$ as $\rho(x)\to0$. Moreover,  for some $c_{18}>1$, there holds
$$\frac1{c_{18}}\rho^s\leq  \xi_1\leq c_{18}\rho^s\quad {\rm in}\ \, \Omega. $$
Using $\xi_1 $ as a test function, we have that
  \begin{eqnarray*}
  \int_{\Omega} u^p\xi_1\,  dx+\int_{\Omega}U_0\xi_1\,   dx&=&\lambda_1\int_{\Omega}w \xi_1\,dx
  \\&\leq &\lambda_1 \Big(\int_{\Omega}w^p \xi_1\,dx\Big)^{\frac1p} \Big(\int_{\Omega}  \xi_1\,dx\Big)^{1-\frac1p},
    \end{eqnarray*}
which implies
\begin{eqnarray*}
  \int_{\Omega}w^p\rho^s\,dx \leq c_{17},
    \end{eqnarray*}
where $c_{17}>0$ depends on $h$. \hfill$\Box$\medskip

 \noindent{\bf Proof of Proposition \ref{pr 2.1}. }
 Suppose by contradiction that there exists a sequence of points $\{x_k\}\subset B_{r_0}\setminus\{0\}$ and a sequence of solutions $u_k$ of (\ref{eq 1.1}) (we can take the same function if $\displaystyle \limsup_{|x|\to0^+} |x|^{N-2s}u(x)=+\infty$) such that $|x_k|\to0^+$ as $k\to+\infty$ and
$$
|x_k|^{N-2s}u_k(x_k) \to+\infty\quad{\rm as}\ \, k\to+\infty.
$$
We can choose  $x_k$ again such that
\begin{equation}\label{ll 3.1}
|x_k|^{N-2s}u_k(x_k)=\max_{x\in\Omega\setminus B_{|x_k|}} |x_k|^{N-2s}u_k(x)\to+\infty\quad{\rm as}\ \, k\to+\infty
\end{equation}
by the fact that the mapping
$r\mapsto \max_{x\in\Omega\setminus B_r} \sigma(x)^{N-2s}u_k(x) $
is nondecreasing.

 We denote 
$$\phi_k(x):= \big(\frac{|x_k|}{2}-|x-x_k| \big)  ^{N-2s}  u_0(x)\quad{\rm for}\ \, |x-x_k|\leq \frac{|x_k|}{2}.$$
Let $\bar x_k$ be the maximum point of $\phi_k$ in $B_{\frac{|x_k|}{2}}(x_k)$, 
that is,
$$\phi_k(\bar x_k)=\max_{|x-x_k|\leq \frac{|x_k|}{2}} \phi_k(x).$$
Set 
$$\nu_k=\frac12\big(\frac{|x_k|}{2}-|\bar x_k-x_k|\big),$$
then $0<2\nu_k<\frac{|x_k|}{2}$ and
$$  \frac{|x_k|}{2}-|x-x_k|\geq \nu_k\ \ {\rm for}\ \, |x-\bar x_k|\leq \nu_k.$$
By the definition of $\phi_k$,  for any $|x-\bar x_k|\leq \nu_k$,
$$(2\nu_k)^{N-2s}u_0(\bar x_k)=\phi_k(\bar x_k)\geq \phi_k(x_k)\geq \nu_k^{N-2s} u_0(x), $$ 
which implies that 
$$2 ^{N-2s}u_0(\bar x_k)\geq u_0(x)\quad {\rm for\ any}\ |x-\bar x_k|\leq \nu_k.$$
Moreover, we see that  
\begin{eqnarray*} 
|\bar x_k|^{N-2s} u_0(\bar x_k)\geq 
(2\nu_k)^{N-2s} u_0(\bar x_k)&=&\phi_k(\bar x_k)\geq \phi_k(x_k)\\[1.5mm]&\geq & (\frac{|x_k|}{2})^{N-2s}u_0(x_k)\to+\infty\quad {\rm as}\ \,  k\to+\infty
\end{eqnarray*} 
by the fact that $|\bar x_k|\geq \frac{|x_k|}{2}\geq 2\nu_k$.
Denote
$$W_k(y)=\frac{1}{u_k(\bar x_k)}u_k\big( u_k(\bar x_k) ^{-\frac{1}{N-2s}} y-\bar x_k \big),\quad\forall\, y\in \Omega_k\setminus \{X_k\},$$
where $$\Omega_k:=\Big\{y\in \R^N: \,   u_k(\bar x_k)^{-\frac{1}{N-2s}} y-\bar x_k\in\Omega\Big\}$$ and
$$X_k= u_k(\bar x_k)^{\frac{1}{N-2s}}   \bar x_k.$$

Note that
$$|X_k|=\big(u_k(\bar x_k) |\bar x_k|^{N-2s} \big)^{\frac{1}{N-2s}}\to+\infty\quad {\rm as}\ \,  k\to+\infty.$$
Thus,   we have that for $y\in \Omega_k\setminus \{X_k\}$
\begin{eqnarray*}
(-\Delta)^s W_k(y)&=&\frac{1}{u_k^{p^*}(\bar x_k)}(-\Delta)^s u_k\big(  u_k(\bar x_k) ^{-\frac{1}{N-2s}} y-\bar x_k\big)
\\[1.5mm]&=&\frac{1}{u_k^{p^*}(\bar x_k)}  u_k^{p^*}\big(  u_k(\bar x_k) ^{-\frac{1}{N-2s}} y-\bar x_k\big)
\\[1.5mm]&=&   W_k^{p^*}(y),
\end{eqnarray*}
that is
\begin{equation}\label{eq lim}
(-\Delta)^s W_k(y)=  W_k^{p^*}(y)\quad {\rm for}\ \ x\in \Omega_k\setminus \{X_k\}.
\end{equation}

We claim that there is $c_{18}>0$ independent of $k$ such that
$$\|W_k\|_{L^1_s(\R^N)} \leq c_{18}$$
and for any $\epsilon>0$, there exists $k_1>0$ and $R>0$ such that
\begin{equation}\label{l1 es}
\int_{\R^N\setminus B_R(0)} W_k(y)(1+|y|)^{-N-2s} dy \leq \epsilon.
\end{equation}

In fact,  since   $|\bar x_k|\to0$, we see that
\begin{eqnarray*}
0&\leq&  \int_{\R^N\setminus \Omega_k} W_k(y)(1+|y|)^{-N-2s} dy
\\[1.5mm]&=&\frac{1}{u_k(\bar x_k)}\int_{\R^N\setminus \Omega_k } u_k\big(  u_k(\bar x_k)  ^{-\frac{1}{N-2s}} y-\bar x_k\big) (1+|y|)^{-N-2s} dy
\\[1.5mm]&=& \frac{1}{u_k^{p^*}(\bar x_k)}\int_{\R^N\setminus \Omega }h(z)  (u_k(\bar x_k)^{-\frac{1}{N-2s}}+|z-\bar x_k|)^{-N-2s}dz
\\[1.5mm]&\leq &\frac{1}{u_k^{p^*}(\bar x_k)} \int_{\R^N  }h(z)  (1+|z-\bar x_k|)^{-N-2s}dz
\\[1.5mm]&\leq & \frac{2}{u_k^{p^*}(\bar x_k)} \int_{\R^N  }h(z)  (1+|z|)^{-N-2s}dz
\to0\quad{\rm as}\ k\to+\infty.
 \end{eqnarray*}
 Taking $r_k= |\bar x_k| u_k(\bar x_k)^{\frac{1}{N-2s}}\to+\infty$ as $k\to+\infty$, we obtain that
 \begin{eqnarray*}
0&\leq&  \int_{B_{r_k}(X_k)} W_k(y)(1+|y|)^{-N-2s} dy\\[1.5mm]&\leq &\frac{r_k^{-N-2s}}{u_k(\bar x_k)}\int_{B_{r_k}(X_k)} u_k\big(u_k(\bar x_k)^{-\frac{1}{N-2s}} y -\bar x_k\big) dy
\\[1.5mm]&=& \frac{r_k^{-N-2s}u_k(\bar x_k)^{\frac{N}{N-2s}} }{u_k(\bar x_k)} \int_{  B_{\frac12 |\bar x_k| } } u_k (z)   dz
\\[1.5mm]&\leq & r_k^{-N-2s}u_k(\bar x_k)^{\frac{2s}{N-2s}}  \Big(\int_{  B_{\frac12 |\bar x_k| } } u_k^{p^*} (z)   dz \Big)^{\frac1{p^*}}  \Big(\int_{  B_{\frac12 |\bar x_k| }}dx\Big)^{1-\frac1{p^*}}
\\[1.5mm]&=& c_{19} r_k^{-N-2s}u_k(\bar x_k)^{\frac{2s}{N-2s}} |\bar x_k|^{2s} \Big(\int_{  B_{\frac12 |\bar x_k| } } u_k^{p^*} (z)   dz \Big)^{\frac1{p^*}}
\\[1.5mm]&\leq & c_{19}  r_k^{-N}
  \to 0\quad{\rm as}\ k\to+\infty
 \end{eqnarray*}
 by Lemma \ref{lm 3.1-1}.

 Moreover, $W_k(y)\leq 1$ in $\Omega_k\setminus B_{r_k}$, then
 \begin{eqnarray*}
   \int_{\Omega_k\setminus B_{r_k}(X_k)} W_k(y)(1+|y|)^{-N-2s} dy&\leq & \int_{\Omega_k\setminus B_{r_k}(X_k)}(1+|y|)^{-N-2s} dy
\\&\leq &  \int_{  \R^N} (1+|y|)^{-N-2s} dy
 \end{eqnarray*}
 and
 \begin{eqnarray*}
  \int_{\big(\Omega_k\setminus B_{r_k}(X_k) \big)\setminus B_R} W_k(y)(1+|y|)^{-N-2s} dy \leq   \int_{\R^N\setminus B_{R}}(1+|y|)^{-N-2s} dy
 \leq    c_{20} R^{-2s},
 \end{eqnarray*}
 where $c_{20}>0$ is independent of $k$.
 Then (\ref{l1 es}) holds true and the claim is proved. \medskip

Note that  $0<W_k\leq 2^{N-2s}$ in $B_{\tilde r_k}$,
where
$$\tilde r_k=\nu_km_k^{\frac{1}{N-2s}}=(\nu_k^{\frac{1}{N-2s}}m_k)^{N-2s}\to+\infty\quad{\rm as}\ \, k\to+\infty,  $$
then for any $R>0$, there exists $k_R$ for any $k\geq k_R$
\begin{eqnarray*}
\|W_k\|_{C^{2s+\alpha}(B_{_R })} &\leq& c_{20}\Big( \|W_k\|_{L^1_s(\R^N)}+\|W_k\|_{L^\infty(B_{_{2R}})}+ \|W_k^p\|_{L^\infty(B_{_{2R}})}\Big)\\[2mm]&\leq& c_{20}\Big( \|W_k\|_{L^1_s(\R^N)}+2^{N-2s}\Big),
\end{eqnarray*}
where $\alpha\in(0,s)$ and $c_{20}>0$.

Since  $0<W_k\leq c_0$ in $B_{r_k}$, so for any $R>0$, there exists $k_R$ for $k\geq k_R$
$$\|W_k\|_{C^{2s+\alpha}(B_{\frac{r_k}{2}})} \leq \|W_k\|_{L^1_s(\R^N)}+\|W_k^p\|_{L^\infty(B_{r_k})},$$
where $\alpha\in(0,s)$.

 Since $R$ is arbitrary,  up to subsequence, there exists a nonnegative function $ W_\infty\in L^\infty(\R^N)$ such that as $k\to+\infty$
 \begin{equation}\label{cve1}
 \arraycolsep=1pt
\begin{array}{lll}
 W_k\to W_\infty\quad\ {\rm in}\ \ C^{2s+\alpha'}_{loc}(\R^N),\\[2mm]
 W_k\to W_\infty\quad\ {\rm in}\ \ L^1_s(\R^N),\\[2mm]
 0<W_\infty  \leq   2 ^{N-2s}
 \end{array}
 \end{equation}
 for some $\alpha'\in(0,\alpha)$.

For any $x\in \R^N$, take $R> 4|x|$ and $R\leq r_k$ for $k$ large enough,  note that
\begin{eqnarray*}
\frac1{c_{N,s}}(-\Delta)^s W_k(x)&=&{\rm p.v.}\int_{B_R} \frac{W_k(x)-W_k(y)}{|x-y|^{N+2s}} dy
+W_k(x) \int_{\R^N\setminus B_R} \frac{1}{|x-y|^{N+2s}} dy
\\&&\ \ -\int_{\R^N\setminus B_R} \frac{W_k(y)}{|x-y|^{N+2s}} dy
\\[1.5mm]&=:&E_{1,k}(x)+E_{2,k}(x)-E_{3,k}(x), \ \, {\rm respectively.}
\end{eqnarray*}
For any $\epsilon>0$,  by (\ref{cve1}), we have that
$$\Big|E_{1,k}(x)-{\rm p.v.}\int_{B_R} \frac{w_\infty(x)-w_\infty(y)}{|x-y|^{N+2s}} dy\Big|<\frac\epsilon3$$
and there exists $R_1>0$ such that for $R>R_1$
$$\Big| E_{2,k}(x)-w_\infty(x) \int_{\R^N\setminus B_R} \frac{1}{|x-y|^{N+2s}} dy\Big|<\frac\epsilon3.$$
Furthermore,  by (\ref{l1 es}) we have that
\begin{eqnarray*}
0<E_{3,k}(x)  &\leq &  \int_{\R^N\setminus B_R} W_k(y)(1+|y|)^{-N-2s} dy<\frac\epsilon3
\end{eqnarray*}
for $k>0$ and $R$ large enough.

Note that  $0<W_k\leq 2^{N-2s}$ in $B_{\tilde r_k}$,
where 
$$\tilde r_k=\nu_k u_k(\bar x_k)^{\frac{1}{N-2s}}=(\nu_k^{N-2s}u_k(\bar x_k)^{\frac{1}{N-2s}}\to+\infty\quad{\rm as}\ \, k\to+\infty,  $$
Then for any $R>0$, there exists $k_R$ for any $k\geq k_R$
\begin{equation}\label{bset re}
\|W_k\|_{C^{2s+\alpha}(B_{_R })} \leq c_{20}\Big( \|W_k\|_{L^1_s(\R^N)}+\|W_k^{p^*}\|_{L^\infty(B_{r_k})}\Big)\leq c_{20}\Big( \|W_k\|_{L^1_s(\R^N)}+2^{N-2s}\Big),
\end{equation}
where $\alpha\in(0,s)$ and $c_{20}>0$.
Therefore,  we conclude that
$$\lim_{k\to+\infty}(-\Delta)^s W_k(x)=(-\Delta)^s W_\infty(x)$$
and $W_\infty$ is a classical solution of
\begin{equation}\label{eq lim1}
(-\Delta)^s W_\infty=W_\infty^{p^*}\quad {\rm in}\ \,\R^N
\end{equation}
satisfying
$$0\leq W_\infty\leq 2^{N-2s}.$$

Since $W_\infty(0)=1$, then $w_\infty>0$ in $\R^N$,  thanks to the nonnegative property of $W_\infty$.
By   \cite[Theorem 3 ]{CLO} or  \cite[Theorem 4.5]{CLO1}, problem (\ref{eq lim1})
has no bounded positive solution since $p^*\in(1,\frac{N+2s}{N-2s})$.
This completes the proof. \hfill$\Box$\medskip

The upper bound in Proposition \ref{pr 2.1} is to obtain the Harnack inequality for singular solution
of (\ref{eq 1.1}).

\begin{proposition}\label{pr Harnack}
 Let $u$ be a nonnegative solution of (\ref{eq 1.1}). Then there exists $C_0>0$ independent of $u$ such that for all $r\in(0,r_0]$
$$
 \sup_{x\in B_{2r}\setminus B_r} u(x)\leq C_0\Big(\inf_{x\in B_{2r}\setminus B_r} u(x)+\|u\|_{L^1_s(\R^N)}\Big).
$$
 If additionally $u$ is singular at the origin, then for all $r\in(0,r_0]$
 \begin{equation}\label{harn in}
 \sup_{x\in B_{2r}\setminus B_r} u(x)\leq C_0 \inf_{x\in B_{2r}\setminus B_r} u(x).
 \end{equation}
 \end{proposition}
{\bf Proof. }    Note that  $B_{r_0}\subset \Omega$ for some $r_0>0$ and without loss of the generality, we set $r_0=1$.   Now for fixed $x_0\in \R^N$ verifying $|x_0|=\frac12$, from there exists $C>0$  independent of $u_0$  such that for any $t\in(0,\frac14]$
 $$\sup_{x\in B_t(x_0)} u_0(x)<C.$$

 Let
 $$w(x)=u_0(x)(1-\eta_0(4x))-h(x) \ \ \ {\rm in} \ \ \R^N,$$
 where $\eta_0$ be a smooth function such that
  $\eta_0=1$ in $B_1(0)$ and $\eta_0=0$ in $\R^N\setminus B_2$,   we recall that $h=0$ in $B_{\frac34}$.  Then
 $$(-\Delta)^s w=|x|^\theta u_0^{p^*-1}w+(-\Delta)^s h+(-\Delta)^s (u_0\eta_0(4x))\quad {\rm for}\ \ x\in B_t(x_0),  $$
  where $(-\Delta)^s h>0$ in $B_1$.
 Note that
  $0<u_0^{p^*-1}\leq c_{15}$ for some $c_{15}>0$ independent of $u_0$ and
\begin{eqnarray*}
0&<&(-\Delta)^s h(x)+(-\Delta)^s (u_0(x)\eta_0(4x))
\\[1.5mm]&\leq& c_{22}\Big(\|h\|_{L^1_s(\R^N)}+\|u_0\|_{L^1_s(\R^N)}\Big)
\\[1.5mm]&\leq&  2c_{22}\|u_0\|_{L^1_s(\R^N)}.
\end{eqnarray*}
Then \cite[Theorem 1.1]{TX} (also see \cite[Theorem 11.1]{CS1}) implies that
$$\sup_{x\in B_t(x_0)} u_0(x)\leq C_1\Big(\inf_{x\in B_t(x_0)} u_0(x)+\|u_0\|_{L^1_s(\R^N )}\Big),$$
which infers
 \begin{equation}\label{harn in00}
 \sup_{x\in B_{2r_0}\setminus B_{r_0}} u_0(x)\leq C_1\Big(\inf_{x\in B_{2r_0}\setminus B_{r_0}} u_0(x)+\|u_0\|_{L^1_s(\R^N)}\Big)
\end{equation}
 by finite covering argument, the scaling property and the upper bound of $u_0$.

 Now we do the scaling:
 $$u_r(x)=r^{N-2s} u_0(rx) $$
 for $r\in(0,\frac14]$.
 Then $u_r$ also verifies (\ref{eq 1.1}) and from Proposition \ref{pr 2.1}, we have
 $$r^{2s} u_r(x)^{p^*-1}\leq C \quad {\rm for}\ \  r<|x|<2r,$$
 where $C$ is dependent of $t$.

 It follows by (\ref{harn in00}) that
 $$
 \sup_{x\in B_{2r }\setminus B_{r }} u_0(x)\leq C_0\Big(\inf_{x\in B_{2r}\setminus B_{r }} u_0(x)+\|u_0\|_{L^1_s(\R^N)}\Big).
$$
 Thanks to
$$\displaystyle \lim_{|x|\to0^+}  u_0(x)=+\infty,$$
we obtain   (\ref{harn in}).  This completes the proof.  \hfill$\Box$

 \subsection{ Lower bound}
 \begin{proposition}\label{pr 2.2}
 Let $u$ be a nonnegative solution of (\ref{eq 1.1}) with non-removable singularity
   $$\lim_{|x|\to0^+} u(x)=+\infty.$$
   Then
  for any $\tau\in(2s-N,0)$
 \begin{equation}\label{ll 3.0}
 \lim_{|x|\to0^+}u(x) |x|^{-\tau}=+\infty.
 \end{equation}

 \end{proposition}
 \noindent{\bf Proof. }
By contradiction, we suppose that  (\ref{eq 1.1}) has a   solution   $u_0$  
such that 
$$\liminf_{|x|\to0^+}u_0(x)|x|^{-\tau}<+\infty$$
for  $\tau\in(2s-N,0)$. 
This gives 
$$\limsup_{|x|\to0^+}u_0(x)|x|^{-\tau}<+\infty$$
by Harnack inequality (\ref{harn in}).

Note that for some $r>0$ and $d_0>0$
$$  u_0^{p^*}(x)\leq d_0 |x|^{\tau p^* } \quad  {\rm for}\ x\in  B_r\setminus\{0\},$$
 where
 $$\tau p^*+2s<N. $$

Let $\tau_0=\tau>2s-N$ and
$$\tau_1:=p^*\tau_0+ 2s,$$
then
$$\tau_1-\tau_0=\frac{2s}{N-2s}\big(\tau-(N-2s)\big)>0.$$
  If $\tau_1>0$, by Lemma \ref{teo 2.3} part $(ii)$, we know that
 $$u_0(x)\leq c_{23}\quad {\rm for}\ x\in B_r\setminus\{0\} $$
 which ends the proof since
 $$\lim_{|x|\to+\infty} u_0(x)=+\infty.$$
If $\tau_1\in\big(2s-N, 0\big]$,  by Lemma \ref{teo 2.3} part $(i)$, we have that
$$u_0(x)\leq d_1|x|^{\tau_1}\quad  {\rm for}\ x\in  B_{r_0}(0)\setminus\{0\}.$$
Iteratively, we recall that
$$\tau_j:=p\tau_{j-1}+\theta+2s,\quad j=1,2,\cdots .$$
Note that
$$\tau_1-\tau_0=(p-1)\tau_0+\theta+2s >0$$
thanks to $\tau_0>-\frac{2s+\theta}{p-1}$.\smallskip

If $\tau_{j}p   +\theta +2s \leq \tau_+(s,\mu)$ the proof is complete, otherwise,
it follows by Lemma \ref{teo 2.3}  that
$$u_0(x)\geq d_{j+1} |x|^{\tau_{j+1}},$$
where
$$\tau_{j+1}=p^*\tau_j+2s<\tau_j.$$

 We claim that  $\{\tau_j\}_j$ is a increasing sequence and
 there exists $j_0\in\N $  such that
 \begin{equation}\label{2.3-bp}
 \tau_{j_0}\le 0 \quad {\rm and}\quad \tau_{j_0-1}>0.
 \end{equation}

In fact,  for   $\tau_0>2s-N$,
$$
\tau_j-\tau_{j-1} = p^*(\tau_{j-1}-\tau_{j-2})=(p^*)^{j-1} (\tau_1-\tau_0)\to+\infty\ \, {\rm as}\ j\to+\infty,
$$
which implies that the sequence $\{\tau_j\}_j$ is  increasing. \medskip

 Thus there exists $j_0\in\N$ such that
$$\tau_{j_0}\leq 0,\quad \quad \tau_{j_0+1}>0 $$
and
\begin{eqnarray*}
 (-\Delta)^s u_0(x) \leq  c_{23}  d_{j_0}^p  |x|^{\tau_{j_0+1}-2s}  \quad{\rm in}\ \ B_{r_0}(0)\setminus\{0\},
 \end{eqnarray*}
 then $\sup_{x\in\Omega\setminus\{0\}} u_0(x)<+\infty$
 and  which contradicts the fact that $u_0$ has non-removable singularity at the origin.   \hfill$\Box$\bigskip

With the help of Proposition \ref{pr 2.1} and Proposition  \ref{pr 2.2}, we are in a position to show Theorem \ref{teo 3-1} \medskip

 \noindent{\bf Proof of Theorem \ref{teo 3-1}. }
 The upper bound and lower bound follow Proposition \ref{pr 2.1} and Proposition  \ref{pr 2.2}
 respectively.
   We complete the proof.\hfill$\Box$

 \setcounter{equation}{0}
  \section{Improved singularity}

  \subsection{Important estimates}
For $R>0$,  let $u_0$ be a positive solution of
 \begin{equation}\label{eq 5.1}
\left\{
\begin{array}{lll}
(-\Delta)^s u=   u^{p^*}  \quad &{\rm in}\ \, B_{R} \setminus\{0\},
\\[1.5mm]
\qquad \ \ \, u=0&{\rm in}\ \,  \R^N\setminus B_{R},\\[1.5mm]
\displaystyle \lim_{|x|\to0^+} u=+\infty.
\end{array}\right.
 \end{equation}
In this section, we will improve the isolated singularity of $u_0$ at the origin.
 \begin{theorem}\label{teo 5.1-1}
For any $R>0$,  any positive solution $u_0$ of (\ref{eq 5.1}) is radially symmetric, strictly decreasing with respect to
 $|x|$.
 \end{theorem}
\noindent  The proof of Theorem \ref{teo 5.1-1} is addressed in  {\it Appendix B}. \medskip

 We consider the function
\begin{equation}\label{def w0}
w_0=v_{-m_0} u_0\quad {\rm in}\ \, \R^N,
\end{equation}
where   $u_0$ is a positive solution of (\ref{eq 5.1}),  $v_{-m_0}$ is  be a smooth,  radially symmetric    function non-increasing with respect to $|x|$   satisfying (\ref{cl-6-1}) 
with
$$m_0=-\frac{N-2s}{2s}<0.$$
%It is worth noting that $v_{-m_0}$ is radially symmetric and decreasing with to $|x|$.

Direct computation shows that   $w_0$ verifies
\begin{equation}\label{eq 4.1-1}
 \left\{ \arraycolsep=1pt
\begin{array}{lll}
 \displaystyle   (-\Delta)^s    w_0= v_{-m_0} u_0^{p^*}-L_s u_0+Q_1 u_0   \quad
 &{\rm in}\quad B_R\setminus\{0\},\\[2mm]
 \phantom{ (-\Delta)^s     }
 \displaystyle    w_0 =0\quad
 &{\rm in}\quad \R^N\setminus B_R,
\end{array}\right.
\end{equation}
where
$$Q_1(x)= (-\Delta)^s v_{-m_0}(x)$$ and
\begin{eqnarray*}
L_su_0(x)  &=&
 c_{N,s} \int_{\R^N} \frac{\big(u_0(x)-u_0(y)\big)\big(v_{-m_0}(x)- v_{-m_0}(y) \big)}{|x-y|^{N+2s}} dy.
\end{eqnarray*}

\begin{lemma}\label{lm est 12}
Let   $u_0$ be a nonnegative classical solution of (\ref{eq 5.1}),
$w_0$ be given in (\ref{def w0})
and
\begin{equation}\label{pt 2}
F(x)=v_{-m_0} u_0^{p^*}-L_s u_0+Q_1 u_0\quad {\rm in}\ \, B_R.
\end{equation}

For any   $\kappa_0\in\big(0, (-\cB_{-m_0})^{\frac1{p^*-1}}\big)$, there exists  $\bar r\in(0,\min\{\frac R2,\frac1{e^2}\}]$ such that
if $$u_0(x)\leq \kappa_0|x|^{2s-N}(-\ln|x|)^{m_0} \quad {\rm in}\ \, 0<|x|<\bar  r,$$
then
\begin{equation}\label{pt 2--1}
 F(x) \leq 0\quad {\rm for}\ \, x\in B_{\bar  r}\setminus\{0\}
  \end{equation}
  and
 $u_0$ is bounded.
\end{lemma}
\noindent{\bf Proof. }   We see that the functions $v_{-m_0},\, u_0$
are radially symmetric and decreasing with respect to $|x|$ by the definition of $v_{-m_0}$ and Theorem \ref{teo 5.1-1}. Thus, we have that
 $$\big(u_0(x)-u_0(y)\big)\big(v_{-m_0}(x)- v_{-m_0}(y) \big)>0\quad {\rm if}\ \ |y|\not=|x|$$
and
$$-L_su_0<0.$$

  From the upper bound of $u_0$ in our assumption, we have that
  for $0<|x|<\min\{\frac R2,\frac{1}{e^2}\}$
  $$v_{-m_0}(x)u_0^{p^*}(x)\leq  \kappa_0^{p^*} (-\ln|x|)^{-1}|x|^{-N}.$$
 Using  (\ref{cl-4}),  we have that
$$\Big|  (-\Delta)^s v_{-m_0} -\cB_{-m_0}  |x|^{-2s}(-\ln |x|)^{\frac{N-2s}{2s}-1}\Big|\leq  c_{24}|x|^{-2s}(-\ln |x|)^{\frac{N-2s}{2s}-2},\ \forall\, 0<|x|<r_0.$$
Thus, for any $\epsilon>0$, there exists $\bar r\in(0, \min\{\frac R2, r_0\})$  such that
\begin{equation}\label{pt 1}
 Q_1(x)\leq  (\cB_{-m_0}+\epsilon)  |x|^{-2s}(-\ln |x|)^{\frac{N-2s}{2s}-1}\quad{\rm for}\ \, 0<|x|<\bar r,
\end{equation}
where $\cB_{-m_0}<0$.
 This implies that for $0<|x|<\bar r$
 $$Q_1(x) u_0(x)\leq  \kappa_0 (\cB_{-m_0}+\epsilon) (-\ln |x|)^{-1}|x|^{-N}.$$
 Note that
 $$v_{-m_0} u_0^{p^*} +Q_1 u_0 -\cL_su_0  \leq  \kappa_0(\kappa_0^{p^*-1}+\cB_{-m_0}+\epsilon)(-\ln |x|)^{-1}|x|^{-N},\  \, x\in B_{\bar r}\setminus\{0\}$$
 by an appropriate choice of $\epsilon$ could taking $\kappa_0^{p^*-1}+\cB_{-m_0}+\epsilon<0$ and that $\kappa_0< (-\cB_{-m_0})^{\frac1{p^*-1}}$.
 Thus, we have that $F\leq 0$ in $B_{\bar r}\setminus\{0\}$ and $F_+=\max\{F,0\}$ is bounded in $B_R\setminus\{0\}$, where $F_+=\max\{F,0\}. $

  Denote
$\mathbb{G}_s[f]$ the Green operator of $f\in L^1(B_R)$ by
$$\mathbb{G}_s[f](x)= \int_{B_R}G_s(x,y) f(y)dy\quad{\rm for}\ \, x\in B_R,$$
where $G_s(\cdot,\cdot)$ is the Green kernel of $(-\Delta)^s$ subject to
the zero Dirichlet condition in $\R^N\setminus B_R$.    There is some $c_{24}>0$ independent of $R$ such that  $G_s(x,y)\leq c_{24}|x-y|^{2s-N}$ for $x\not=y$.  Here we refer  to \cite{CS} for the properties of Green kernel.  Thus  $\mathbb{G}_s[F_+]$ is bounded.      

 From  Lemma \ref{cr hp} we have that
\begin{eqnarray*}
 0\leq w_0(x)  \leq  \mathbb{G}_s[F_+](x).
\end{eqnarray*}
Consequently,
 $$ u_0(x)\leq  (-\ln|x|)^{-m_0}\quad {\rm in}\ \, B_{\bar r}\setminus\{0\}.$$
For some $\tau'\in(2s-N,0)$ and $c_{24}>0$, we have:
 \begin{equation}\label{bound e-1}
 u_0(x)\leq  c_{24} |x|^{\tau'}\quad {\rm in}\ \, B_{\bar r}\setminus\{0\}.
 \end{equation}

  Now we prove that  $u_0$ is bounded.  If not, we can assume that 
$$\limsup_{|x|\to0^+} u_0(x)=+\infty,$$
which implies by the Harnack inequality  Proposition \ref{pr Harnack} that  
$$\lim_{|x|\to0^+} u_0(x)=+\infty$$
From Theorem \ref{teo 3-1},
we have that
$$\lim_{|x|\to0^+} u_0(x)  |x|^{-\tau'}=+\infty,$$
which contradicts (\ref{bound e-1}). So $u_0$ is bounded. 
 \hfill$\Box$\medskip

To improve the upper bound,  we need to consider the Liouville type theorem for the fractional Poisson problem  with a weak Hardy potential  
\begin{equation}\label{eq 1.1-qh}
\left\{
\begin{array}{lll}
(-\Delta)^s u= \frac{\nu}{|x|^{2s}(-\ln|x|)} u+   f \quad &{\rm in}\ \, B_{R_1} \setminus\{0\},
\\[2mm]
\qquad \ \   u\geq 0&{\rm in}\ \,  \R^N\setminus B_{R_1},
\end{array}\right.
 \end{equation}
where $R_1\in(0,\frac1{e}]$ and $f:B_{R_1}\setminus\{0\}\to \R$ is H\"older continuous locally in $\bar B_{R_1}\setminus\{0\}$.

 \begin{lemma}\label{teo 2.1-qh} Let $\nu>0$, and $f$ be a nonnegative function such that $f\in C^\beta_{loc}(\bar B_{R_1} \setminus \{0\})$ for some $\beta\in(0,1)$.
The homogeneous problem (\ref{eq 1.1-qh}) has no positive solution if
\begin{equation}\label{as-qh}
\liminf_{|x|\to0^+} f(x) |x|^{N}(-\ln|x|)^{1+\frac{\nu}{-\cC_s'(0)}} >0.
\end{equation}
\end{lemma}
{\bf Proof. }  By contradiction, we assume that problem (\ref{eq 1.1-qh})
has a positive solution $u_0$.  

From the assumption  (\ref{as-qh}),   we take
\begin{equation}\label{eq-qh-0}
m_1=\frac{\nu}{\cC_s'(0)}\in(-\infty,0).
\end{equation}
 Let us set:
$$ \liminf_{|x|\to0^+} f(x) |x|^{N}(-\ln|x|)^{1-m_1}=2\bar k_0,$$
then  there exists
  $r_1\in(0, \min\{r_0,R_0\})$ such that
$$f(x)\geq  \bar k_0  |x|^{-N}(-\ln|x|)^{m_1-1} \quad {\rm for}\ \ 0<|x|<r_1, $$
where $\cB_{m_1}>0$ for $m_1<0$.
Let $W_0$ be the solution of
\begin{equation}\label{eq 1.1-qh-0}
\left\{
\begin{array}{lll}
(-\Delta)^s u=  \bar k_0   |x|^{-N}(-\ln|x|)^{m_1-1}\chi_{_{B_{r_1}}} \quad &{\rm in}\ \, B_{R_1} \setminus\{0\},
\\[2mm]
\qquad \ \   u=0&{\rm in}\ \,  \R^N\setminus B_{R_1},
\end{array}\right.
 \end{equation}
 where $\chi_{_{B_{r_1}}}$ is the characterized function of $B_{r_1}$.
By  direct comparison with $w_{m_1}(R_0^{-1}\cdot)$ and Lemma \ref{lm 2.1-d} with
$m=m_1$ and $\cD_{m_1}<0$, by re-choice of $r_1$ if necessary, we have that
%there exists $\epsilon_0\in(0, \frac{\bar k_0}{4})$ such that
$$W_0(x)\geq \frac{\bar k_0 }{\cB_{m_1}} |x|^{2s-N}(-\ln|x|)^{m_1}\quad {\rm for}\  \,  0<|x|<r_1.$$
Then Lemma \ref{cr hp}  implies that
$$u_0(x)\geq W_0(x)\geq  \frac{\bar k_0 }{\cB_{m_1}}   |x|^{2s-N}(-\ln|x|)^{m_1}\quad {\rm for}\ \ 0<|x|<r_1  $$
and now let
$$k_0=\frac{\bar k_0 }{\cB_{m_1}}, $$
then
\begin{eqnarray*}
\frac{\nu}{|x|^{2s}(-\ln|x|)} u_0(x) &\geq&   k_0 \nu   |x|^{-N} (-\ln|x|)^{m_1-1}\quad {\rm for}\  \,  0<|x|<r_1.
 \end{eqnarray*}

So we have that
$$(-\Delta)^s u_0\geq  ( \bar k_0   +k_0\nu )   |x|^{-N} (-\ln|x|)^{m_1-1}\quad {\rm for}\ 0<|x|<r_1.$$

 Then by Lemma \ref{cr hp} again,  we have that
 $$u_0\geq  k_1  W_0 \quad {\rm for}\ \ 0<|x|<R_0, $$
 where
 $$k_1= \frac{\bar k_0   +k_0\nu }{\cB_{m_1}}=k_0+ \frac{ \nu }{\cB_{m_1}}k_0$$
and then for $0<|x|<r_1$
\begin{eqnarray*}
\frac{\nu}{|x|^{2s}(-\ln|x|)} u_0(x) &\geq& k_1k_0    |x|^{-N} (-\ln|x|)^{m-1}.
 \end{eqnarray*}
By repeating the above procedure,  we have that
$$u_0\geq k_n W_0 \quad {\rm for}\ \ 0<|x|<r_1, $$
where
$$k_n=  k_0+ \frac{ \nu }{\cB_{m_1}}k_{n-1}. $$
Note that
\begin{eqnarray*}
k_n- k_{n-1}=\frac{\nu   }{\cB_{m_1}}  (k_{n-1}-k_{n-2}) =  (\frac{\nu   }{\cB_{m_1}})^{n}k_0= k_0,
  \end{eqnarray*}
where $\frac{\nu   }{\cB_{m_1}}=1$ by the choice of $ m_1$ in (\ref{eq-qh-0})
and $\cB_{m_1}=\cC_s'(0)m_1$.
 Then  $$k_n \to+\infty\quad {\rm as}\ \, n\to+\infty,$$
which implies that   $u_0$ blows up in $B_{r_1}$ and we obtain a contradiction.   \hfill$\Box$

 \subsection{Isolated singularity of (\ref{eq 1.1}) }

 \begin{lemma}\label{lm existence 1}
 Let   $h $ be a nonnegative function in $ C^\theta(B_{2R_0})\cap L^1_s(\R^N)$ with $\theta>2s$
 and  $u_0$ be a positive singular solution of  (\ref{eq 1.1}) such that
  \begin{eqnarray*}
\limsup_{|x|\to0^+}u_0(x) (|x|(-\ln|x|)^\frac{1}{2s})^{N-2s}<+\infty.
\end{eqnarray*}

Then for $r_1\in(0, \frac1{e^2}]$,   there exists  $m>0$  such that
 problem (\ref{eq 5.1}) with $r=r_1$ has a positive  singular solution $u_1$ such that
 $$u_0(x)- (-\ln |x|)^{m+1}\leq u_1(x)\leq u_0(x)\quad {\rm in}\ \, B_{r_1}\setminus\{0\}.$$
 \end{lemma}
\noindent{\bf Proof. } From the upper bound assumption, there exists $c_{25}>0$ such that
$$u_0(x)\leq c_{25}|x|^{2s-N}(-\ln|x|)^{m_0}\quad {\rm for} \ \, 0<|x|<r_1. $$

 Let
 $$v_0=u_0\quad {\rm in}\ \, \R^N .$$
We denote $v_n$   the solution of
  \begin{equation}\label{eq 6.1-n}
\left\{
\begin{array}{lll}
(-\Delta)^s v_n=   v_{n-1}^{p^*}  \quad\,  {\rm in}\ \, B_{r_1} \setminus\{0\},
\\[1.5mm]
\qquad \ \ \, v_n=0 \quad\ \ \quad {\rm in}\ \,  \R^N\setminus B_{r_1},\\[1.5mm]
\displaystyle \lim_{|x|\to0^+} v_n(x)|x|^{N-2s}=0.
\end{array}\right.
 \end{equation}
Obviously, we have that
$$(-\Delta)^s (v_1-v_0)\leq 0$$
 and the comparison principle implies that $0<v_1\leq v_0$ in $\R^N\setminus\{0\}$.
 Inductively,  the mapping
 $$n\in\N\mapsto v_n$$
 is decreasing.
 Therefore, by stability results \cite[Theorem 2.4]{CFQ} (also see   \cite[Lemma 4.3]{CS2} for bounded sequence) and the regularity result \cite[Theorem 2.1]{CFQ}, the limit of $\{v_n\}$, $v_\infty$ exists,
 $$\lim_{n\to+\infty} v_n(x)=v_\infty(x)\quad {\rm for\ any}\ \, x\in\R^N\setminus\{0\}, $$
 then $v_\infty$ is a solution of  (\ref{eq 5.1})
 and it verifies that
 $$0\leq v_\infty\leq v_0\quad {\rm in}\ \, \R^N\setminus\{0\}.$$

Let $\nu_1=v_0-v_1$, then
  \begin{eqnarray*}
  (-\Delta)^s \nu_1\leq 0 \quad\  {\rm in}\ \, B_{r_1} \setminus\{0\}
  \end{eqnarray*}
and
$$\nu_1=u_0 \quad\  {\rm in}\ \,  \R^N\setminus B_{r_1}, \qquad
 \lim_{|x|\to0^+} \nu_1(x)|x|^{N-2s}=0.
$$
 Then there exists $M_0>0$ depending on $r_1$ such that
 $$0<\nu_1(x)\leq M_0.$$
For $n=2,3,\cdots$,  denote
$$\nu_n=v_{n-1}-v_n,$$
then
 \begin{eqnarray*}
 (-\Delta)^s \nu_n(x) &=& v_{n-1}^{p^*}(x)-v_{n}^{p^*}(x)
  \\[1.5mm]&\leq& v_0^{p^*-1}(x)\nu_{n-1}(x)
 \leq  c_{25}^{p^*-1} |x|^{-2s}(-\ln|x|)^{-1} \nu_{n-1}(x).
 \end{eqnarray*}
 In particular, we have that
  \begin{eqnarray*}
 (-\Delta)^s \nu_2(x) \leq c_{25}^{p^*-1} M_0 |x|^{-2s}(-\ln|x|)^{-1}  .
 \end{eqnarray*}
Let  $r_1\leq \frac1{e^2}$ and then for $-\ln r_1\geq 2$ and   $\epsilon_0>0$,  there exists $m>0$  such that
\begin{eqnarray*}
 (-\Delta)^s \nu_2(x) \leq \epsilon_0 |x|^{-2s}(-\ln|x|)^{m}\quad {\rm for}\ \, 0<|x|<r_1.
 \end{eqnarray*}
  Here $m>0$ is chosen such that
  $$\epsilon_0\geq (-\ln r_1)^{-m+1} c_{25}^{p^*-1} M_0. $$

Note that   for $0<|x|<r_1$
  \begin{eqnarray*}
  (-\Delta)^\s v_{m+1}(x) &=&  -\cB_{m+1}  |x|^{-2s}(-\ln |x|)^{m}+O(1)|x|^{-2s}(-\ln|x|)^{m-1}
 \\[1.5mm]& \geq& (-\cB_{m+1}-\epsilon_0)   |x|^{-2s}(-\ln |x|)^{m},
\end{eqnarray*}
where $-\cB_{m+1}>3\epsilon_0>0$ since $m+1>0$.

 Therefore, we have that
  $$0<\nu_2(x)\leq M_1  (-\ln|x|)^{m+1}, $$
  where
  $$M_1=   \frac{ \epsilon_0  }{-\cB_{m+1}-\epsilon_0}<\frac12.$$

 By repeating the above process, we can obtain that
 $$0<\nu_n(x)\leq M_n  (-\ln|x|)^{m+1}, $$
  where
  $$M_n=M_{n-1}  \frac{ \epsilon_0  }{\cB_{m+1}-\epsilon_0} =\Big( \frac{ \epsilon_0  }{\cB_{m+1}-\epsilon_0}\Big)^n.$$
Since
 $$0< \frac{ \epsilon_0  }{\cB_{m+1}-\epsilon_0}<\frac12,$$
 then
\begin{eqnarray*}
v_0\geq v_\infty(x)&\geq& v_0(x)-\sum_{n=1}^{+\infty} \nu_n(x)
\\[1.5mm]&\geq& v_0(x)-  \sum^\infty_{n=1} \Big(\frac{  \epsilon_0 }{ \cB_{ 1}-\epsilon_0}\Big)^n   (-\ln|x|)^{m+1}
\\[1.5mm]&\geq&v_0(x)-   (-\ln|x|)^{m+1}.
\end{eqnarray*}
 Which completes the proof.\hfill$\Box$\medskip

\noindent{\bf Proof of Theorem \ref{teo 5.1}. }   Our proof is divided into two parts. \smallskip

{\it Part I: the Ball $\Omega=B_1$ and $h\equiv 0$.}  Let $u_0$ be a positive solution of (\ref{eq 1.1})
and recall that
$$\cK_s=  \cB_{m_0}^{\frac{N-2s}{2s}}.$$
  Set
\begin{equation}\label{1.2-400}
k= \liminf_{|x|\to0^+}u_0(x) \big(|x|(-\ln|x|)^\frac{1}{2s}\big)^{N-2s}.
 \end{equation}
If $k=0$, from the Harnack inequality, we have that
$$\lim_{|x|\to0^+}u_0(x) \big(|x|(-\ln|x|)^\frac{1}{2s}\big)^{N-2s}=0$$
 and Lemma \ref{lm est 12} implies that $u_0$ is bounded.

So if $u_0$ has a non-removable singularity at the origin, the Harnack inequality implies
$$\lim_{|x|\to0^+}u_0(x)=+\infty,$$
then we have that $k>0$.\smallskip

  Now we claim
 $$  k\leq \cK_s. $$
In fact, if $$k\in(\cK_s,+\infty],$$
 letting
$$\epsilon_0=\frac{k-\cK_s}{2\cK_s},$$
then by Lemma \ref{lm 2.2-lower}, there exists $r_1\in(0,\frac1{e^2})$ such that
$$u_0(x)\geq \cK_s\big(1+\epsilon_0)(|x|(-\ln|x|)^\frac{1}{2s}\big)^{2s-N} \quad {\rm for}\ \ 0<|x|<r_1, $$
then for $0<|x|<r_1$
\begin{eqnarray*}
 u_0^{p^*}(x)&\geq& \cK_s^{p^*} \big(1+\epsilon_0)^{p^*} (|x|(-\ln|x|)^\frac{1}{2s}\big)^{-N}
\\[2mm]&=&  \cB_{m_0}(1+\epsilon_0)^{p^*}(|x|(-\ln|x|)^\frac{1}{2s}\big)^{-N} .
 \end{eqnarray*}

 Therefore,     $u_0$ verifies that
\begin{equation}\label{eq 6.1-dd}
\left\{
\begin{array}{lll}
(-\Delta)^s u_0= \frac{\cB_{m_0}}{|x|^{2s}(-\ln|x|)} u_0+ f    \quad &{\rm in}\ \, B_{r_0} \setminus\{0\},
\\[2mm]
\qquad \ \   u_0\geq 0&{\rm in}\ \,  \R^N\setminus B_{r_0},
\end{array}\right.
\end{equation}
where
\begin{eqnarray*}
f(x) = \Big(u_0^{p^*-1}-\frac{\cB_{m_0}}{|x|^{2s}(-\ln|x|)}\Big) u_0(x)
 \geq  \epsilon_0 |x|^{-N}(-\ln|x|)^{m_0-1}.
\end{eqnarray*}
Then
$$\liminf_{|x|\to0^+} f(x)|x|^{N}(-\ln|x|)^{1-m_0} >0$$
and a contradiction arises  Lemma  \ref {teo 2.1-qh} with $\nu=\cB_{m_0}$, from which problem (\ref{eq 6.1-dd}) has no such positive solution. Therefore we obtain  that
  $  k\leq \cK_s. $

  Set
\begin{equation}\label{1.2-40}
\kappa= \limsup_{|x|\to0^+}u_0(x) \big(|x|(-\ln|x|)^\frac{1}{2s}\big)^{N-2s}.
 \end{equation}
 Let us prove that
 $$  \kappa\geq \cK_s. $$
In fact, if
$$\kappa<\cK_s, $$ letting
$$\epsilon_1=\min\Big\{\frac{\cK_s-\kappa}{2\cK_s},\, \frac14\Big\},$$
then by Lemma \ref{lm 2.2-upper} there exists $r_1\in(0,\frac1{e^2})$ such that
$$u_0(x)\leq \cB_{m_0}\big(1-\epsilon_1)(|x|(-\ln|x|)^\frac{1}{2s}\big)^{2s-N} \quad {\rm for}\ \ 0<|x|<r_1, $$
then for $0<|x|<r_1$
\begin{eqnarray*}
 u_0^{p^*}(x)&\leq& \cB_{m_0}^{p^*}\big(1-\epsilon_1)^{p^*} (|x|(-\ln|x|)^\frac{1}{2s}\big)^{-N}.
 \end{eqnarray*}
Set $$s_0\in\Big(0,\frac{p^*-1}{4}\Big),$$ then for $r_1$ small enough,
$$(-\Delta)^s w_{m_0}(x)\geq    \cB_{m_0}\big(1-\epsilon_1)^{s_0} |x|^{-N} (-\ln|x|)^{-\frac{N}{2s}} \quad {\rm for}\ \ 0<|x|<r_1, $$
since $\cD_{m_0}<0$.
 Then by Lemma \ref{lm 2.2-upper}, there exists $r_2\in(0,r_1]$ such that
 $$u_0\leq \cB_{m_0}(1-\epsilon_1)^{p^*-s_0} (|x|(-\ln|x|)^\frac{1}{2s}\big)^{2s-N}\quad {\rm for}\ \ 0<|x|<r_2, $$
and then
\begin{eqnarray*}
 u_0^{p^*}(x) \leq    \cK_s^{p^*}  (1-\epsilon_1)^{p^*(p^*-s_0)}(|x|(-\ln|x|)^\frac{1}{2s}\big)^{-N} .
 \end{eqnarray*}
By repeating the above procedure, there exists $r_n>0$ such that
$$u_0\leq \cK_s(1-\epsilon_1)^{t_n} (|x|(-\ln|x|)^\frac{1}{2s}\big)^{2s-N}\quad {\rm for}\ \ 0<|x|<r_n, $$
where
\begin{eqnarray*}
t_0=p^* \quad {\rm and}\quad  t_n =p^*(t_{n-1}-s_0).
  \end{eqnarray*}
 Note that
 \begin{eqnarray*}
t_n-t_{n-1}= p^* (t_{n-1}-t_{n-2})&=&(p^*)^{n-1}(t_1-t_0)
\\[1.5mm]&=&(p^*)^{n}(p^*-s_0-1) \to+\infty\quad{\rm as}\ \, n\to+\infty,
  \end{eqnarray*}
 then there exists $n_1\in\N$ such that $\cK_s(1-\epsilon_1)^{t_n}<\kappa$, which contradicts  (\ref{1.2-40}) when $n\geq n_1$.

Now  we conclude that
\begin{eqnarray*}
0 < \liminf_{|x|\to0^+}u(x) (|x|(-\ln|x|)^\frac{1}{2s})^{N-2s}
  \leq  \cK_s   \leq \limsup_{|x|\to0^+}u(x) (|x|(-\ln|x|)^\frac{1}{2s})^{N-2s}\leq +\infty.
\end{eqnarray*}
However, if $\displaystyle \limsup_{|x|\to0^+}u(x) (|x|(-\ln|x|)^\frac{1}{2s})^{N-2s}= +\infty$, the Harnack inequality (\ref{harn in}) implies that
$$\liminf_{|x|\to0^+}u(x) (|x|(-\ln|x|)^\frac{1}{2s})^{N-2s}=+\infty. $$
Therefore, by the Harnack inequality
\begin{eqnarray*}
\frac{\cK_s}{C_0} &\leq& \liminf_{|x|\to0^+}u(x) (|x|(-\ln|x|)^\frac{1}{2s})^{N-2s}
  \leq  \cK_s   
  \\&  \leq &\limsup_{|x|\to0^+}u(x) (|x|(-\ln|x|)^\frac{1}{2s})^{N-2s}\leq \cK_s C_0.
\end{eqnarray*}

By the scaling technique, the classification of singularities for positive solutions of (\ref{eq 1.1})
holds in any ball $B_r$. 

\smallskip

{\it Part II: General domain.} In general bounded domain $\Omega$  
$$B_1\subset \Omega\subset B_{R_0}.$$ 
  Let $u_0$ be a positive solution of (\ref{eq 1.1})
such that
$$\lim_{|x|\to0^+}u_0(x)=+\infty.$$
Note that the argument of the upper bound doesn't use the radially symmetric property of Theorem \ref{teo 5.1-1}, so we have that
$$k_0\leq \limsup_{|x|\to0^+}u_0(x)|x|^{N-2s}(-\ln |x|)^{-m_0}<+\infty.$$

By contradiction, we set
$$\liminf_{|x|\to0^+}u_0(x)|x|^{N-2s}(-\ln |x|)^{-m_0}=0.$$

By Lemma \ref{lm existence 1},
   problem (\ref{eq 5.1}) with $R=R_0$ has a positive  singular solution $u_1$ such that
 $$u_0(x)-c_m(-\ln |x|)^{m+1}\leq u_1(x)\leq u_0(x), \ \ x\in B_{r_0}\setminus\{0\}$$
 for some $m>0$ and $c_m>0$.
Then there holds
$$\limsup_{|x|\to0^+}u_1(x)|x|^{N-2s}(-\ln |x|)^{-m_0}\geq \cK_s$$
and
$$\liminf_{|x|\to0^+}u_1(x)|x|^{N-2s}(-\ln |x|)^{-m_0}=0,$$
which are impossible with our above argument when $\Omega=B_{R_0}$.
 \hfill$\Box$

 \setcounter{equation}{0}
  \section{Existence of singular solutions }

 \subsection{ First singular solution}

 From Proposition \ref{pr 2.1-a} with $m=m_0=-\frac{N-2s}{2s}<0$,
 $$\cB_{m_0}=m_0\cC_s'(0)>0,\qquad \cD_{m_0}=\frac{m_0(m_0-1)}{2}\cC''_s(0)<0,$$
 recalling that 
 $$w_{m_0}(x)=|x|^{2s-N}(-\ln|x|)^{m_0}\quad{\rm for}\  x\in B_{\frac1{e^2}}\setminus\{0\},$$
 there exists $r_1\in(0,r_0)$ such that for $0<|x|<r_1$
\begin{eqnarray*}
  (-\Delta)^\s w_{m_0}(x) &=&  \cB_{m_0}  |x|^{-N}(-\ln |x|)^{m_0-1}+\cD_{m_0}|x|^{-2s}(-\ln|x|)^{m_0-2}
  \\[1.5mm]&& +O(1) |x|^{-N}(-\ln |x|)^{m_0-3}
 \\[1.5mm]& \leq&  \cB_{m_0}  |x|^{-N}(-\ln |x|)^{m_0-1}.
\end{eqnarray*}
Recall that 
 $$\cK_s =   \cB_{m_0} ^{\frac{N-2s}{2s}},$$
 then $w_{m_0}$ verifies that
 \begin{equation}\label{sub 6.1}
  (-\Delta)^\s (\cK_s w_{m_0}) \leq (\cK_s w_{m_0})^{p^*}\quad {\rm in}\ \, B_{r_1}\setminus\{0\}.
 \end{equation}

 \begin{proposition}\label{existence 2}
 Let $N>2s$,   then for $r\in(0,r_1]$
 problem
  \begin{equation}\label{eq 6.1-case1}
\left\{
\begin{array}{lll}
(-\Delta)^s u=   u^{p^*}  \quad &{\rm in}\ \, B_{r} \setminus\{0\},
\\[2mm]
\qquad \ \ \, u=0&{\rm in}\ \,  \R^N\setminus B_{r}
\end{array}\right.
 \end{equation}
 has a positive singular solution $u_r$ verifying for some $m>0$
 $$-\cK_sr^{2s-N}(-\ln r)^{m_0} \leq   u_r(x) - \cK_s|x|^{2s-N}(-\ln|x|)^{m_0}\leq  c_{27} |x|^{2s-N}(-\ln|x|)^{m_0-1},$$
where  $c_{27}>0$ is independent of $r$. 

Moreover, the mapping $r\in(0,r_1] \mapsto u_r$ is increasing.
 \end{proposition}
\noindent{\bf Proof. }
Let
 $$v_0=\cK_s w_{m_0}\quad {\rm in}\ \, \R^N $$
 and  we first show that
 \begin{equation}\label{eq 6.1-case2}
\left\{
\begin{array}{lll}
(-\Delta)^s u=   u^{p^*}  \ \ &{\rm in}\ \, B_{r_1} \setminus\{0\},
\\[1.5mm]
\qquad \ \ \, u=v_0&{\rm in}\ \,  \R^N\setminus B_{r_1}
\end{array}\right.
 \end{equation}
 has a solution $u_{r_1}$ verifying
$$\lim_{|x|\to0^+} u_{r_1}(x)|x|^{N-2s}(-\ln|x|)^{-m_0}=\cK_s.$$

To this end, we set
 $$\epsilon_1\in\Big(0,\min\big\{-\frac{\cD_{m_0}}4,\,  \frac12 \frac{\cB_{m_0-1}-\cB_{m_0}}{\cB_{m_0-1}+1}\big\}\Big),$$
 where $\cB_{m_0-1}-\cB_{m_0}=-\cC_s'(0)>0$.  For $r_1>0$,   we re-choose it smaller if necessary,  there holds that  for $0<|x|<r_1$  
  \begin{eqnarray*}
  (-\Delta)^\s v_0(x) &\leq& \cK_s^{p^*}  |x|^{-N}(-\ln |x|)^{m_0-1}+\cK_s(\cD_{m_0}+\epsilon_1) |x|^{-2s}(-\ln|x|)^{m_0-2}
  \\[1.5mm]&= & v_0^{p^*}(x)+\cK_s(\cD_{m_0}+\epsilon_1) |x|^{-2s}(-\ln|x|)^{m_0-2},
\end{eqnarray*}
where we recall that $\cD_{m_0}=\frac{m_0(m_0-1)}{2}\cC_s''(0)<0$.

Inductively,  we denote $v_n$,   (also use the notation $v_{r,n}$ if $r_1$ is replaced by $r$)    the solution of
$$
\left\{
\begin{array}{lll}
(-\Delta)^s v_n=   v_{n-1}^{p^*}  \quad\ \, {\rm in}\ \, B_{r_1} \setminus\{0\},
\\[1.5mm]
\qquad \ \ \, v_n=v_0 \quad\ \ \quad {\rm in}\ \,  \R^N\setminus B_{r_1},\\[1.5mm]
\displaystyle \lim_{|x|\to0^+} v_n(x)|x|^{N-2s}=0.
\end{array}\right.
$$
Obviously, we have that
$$(-\Delta)^s (v_1-v_0)\geq 0$$
 and the comparison principle implies that $ v_1\geq v_0$ in $\R^N\setminus\{0\}$.
 Inductively,  the mapping
 $$n\in\N\mapsto v_n$$
 is increasing.

Set $\nu_1=v_1-v_0$, then
  \begin{eqnarray*}
  (-\Delta)^s \nu_1\leq  \cK_s(-\cD_{m_0}+\epsilon) |x|^{-2s}(-\ln|x|)^{m_0-2} \quad\  {\rm in}\ \, B_{r_1} \setminus\{0\}
  \end{eqnarray*}
and
$$\nu_1=0 \quad\  {\rm in}\ \,  \R^N\setminus B_{r}, \qquad
 \lim_{|x|\to0^+} \nu_1(x)|x|^{N-2s}=0.
$$
By the comparison Principle, we have that
 $$0<\nu_1(x)\leq Z_1 |x|^{2s-N}(-\ln|x|)^{m_0-1},$$
 where
 $$Z_1=\frac{\cK_s(-\cD_{m_0}+\epsilon_1)}{ \cB_{m_0-1}-\epsilon_1}. $$
For $n=2,3,\cdots$,  denote
$$\nu_n=v_n-v_{n-1}\quad {\rm in}\ \, \R^N\setminus\{0\}.$$
If $p^*\leq 2$
 \begin{eqnarray*}
 (-\Delta)^s \nu_2(x)  =  v_{1}^{p^*}(x)-v_{0}^{p^*}(x)
 &\leq& v_1^{p^*-1}(x)\nu_{1}(x)
 \\[1.5mm]&\leq&   (v_0+\nu_1)^{p^*-1}(x)\nu_{1}(x)
  \\[1.5mm]&\leq&  v_0^{p^*-1}(x)\nu_{1}(x)+ \nu_{1}^{p^*}(x);
 \end{eqnarray*}
 if $p^*>2$
  \begin{eqnarray*}
 (-\Delta)^s \nu_2(x)   &\leq&   (v_0+\nu_1)^{p^*-1}(x)\nu_{1}(x)
  \\[1.5mm]&\leq&  v_0^{p^*-1}(x)\nu_{1}(x)+2^{p^*}v_0^{p^*-2}(x)\nu_{1}^2(x) +2^{p^*}\nu_{1}^{p^*}(x).
 \end{eqnarray*}

In the case $p^*\leq 2$,  for  $0<|x|<r_1$ one has that
 \begin{eqnarray*}
 0<\nu_2(x)&\leq& Z_1 \frac{ \cB_{m_0}}{ \cB_{m_0-1}-\epsilon_1} |x|^{2s-N}(-\ln|x|)^{m_0-1}
 \\[1.5mm]&& + \frac{Z_1^{p^*} }{ \cB_{(m_0-1)p^*+1}-\epsilon_1}  |x|^{2s-N}(-\ln|x|)^{(m_0-1)p^*+1}
\\[1.5mm]&\leq& Z_2|x|^{2s-N}(-\ln|x|)^{m_0-1},
  \end{eqnarray*}
  where
  $$Z_2=Z_1\Big( \frac{ \cB_{m_0}}{ \cB_{m_0-1}-\epsilon_1}+\epsilon_1 \Big)$$
  and we used the fact that
  $$(m_0-1)p^*+1<m_0-1.$$

In the case $p^*>2$,  for  $0<|x|<r_1$
 \begin{eqnarray*}
 0<\nu_2(x)&\leq& Z_1 \frac{ \cB_{m_0}}{\cB_{m_0-1}-\epsilon_1} |x|^{2s-N}(-\ln|x|)^{m_0-1}
 +2^{p^*} Z_1^2|x|^{2s-N}(-\ln|x|)^{m_0-2}
 \\[1.5mm]&&+2^{p^*}\Big(\frac{Z_1^{p^*} }{ \cB_{(m_0-1)p^*+1}-\epsilon_1}\Big) |x|^{2s-N}(-\ln|x|)^{(m_0-1)p^*+1}
\\[1.5mm]&\leq& Z_2|x|^{2s-N}(-\ln|x|)^{m_0-1}.
  \end{eqnarray*}

 By iterating the above process, we can obtain that
 $$0<\nu_n(x)\leq Z_n |x|^{2s-N}(-\ln|x|)^{m_0-1}, $$
  where
  $$Z_n=Z_{n-1} \Big(\frac{ \cB_{m_0}}{ \cB_{m_0-1}-\epsilon_1}+\epsilon_1\Big)=Z_1  \Big(\frac{ \cB_{m_0}}{\cB_{m_0-1}-\epsilon_1}+\epsilon_1\Big)^n.$$
 Note that
\begin{eqnarray*}0&<&\frac{ \cB_{m_0}}{ \cB_{m_0-1}-\epsilon_1}+\epsilon_1
\\[1mm]& =&1- \frac{\cB_{m_0-1}-\cB_{m_0}-\epsilon_1\cB_{m_0-1}+\epsilon_1+\epsilon_1^2}{\cB_{m_0-1}-\epsilon_1} 
\\[1mm]& \leq & 1- \frac{(\cB_{m_0-1}-\cB_{m_0})\big(1-\frac{\cB_{m_0-1}}{2(\cB_{m_0-1}+1)}\big)  + \epsilon_1}{\cB_{m_0-1}-\epsilon_1}
  <  1,
\end{eqnarray*}
where the last inequality holds  by the choice of $\epsilon_1$. Then we have that 
\begin{eqnarray*}
v_0\leq  v_n(x)&\leq& v_0(x)+\sum_{n=1}^{n} \nu_n(x)
\\[1.5mm]&\leq& v_0(x)+Z_1\Big(\sum^n_{n=1} \Big(\frac{ \cB_{m_0}}{ \cB_{m_0-1}-\epsilon_1}+\epsilon_1\Big)^n\Big) |x|^{2s-N}(-\ln|x|)^{m_0-1}
\\[1.5mm]&<& v_0(x)+\bar v(x),
\end{eqnarray*}
where
$$\bar v(x)=Z_1\Big(\sum^{+\infty}_{n=1} \Big(\frac{ \cB_{m_0}}{ \cB_{m_0-1}-\epsilon_1}+\epsilon_1\Big)^n\Big) |x|^{2s-N}(-\ln|x|)^{m_0-1}\quad {\rm for}\ \,  x\in B_{r_1}\setminus\{0\}.$$

  Therefore, by stability results \cite[Theorem 2.4]{CFQ} (also see   \cite[Lemma 4.3]{CS2} for bounded sequence) and the regularity result \cite[Theorem 2.1]{CFQ}, the limit of $\{v_n\}$ exists, denoting  $v_{r_1,\infty}$,
 $$\lim_{n\to+\infty} v_n(x)=v_{r_1,\infty}(x)\quad {\rm in}\ \, \R^N\setminus\{0\}, $$
 then $v_{r_1,\infty}$ is a solution of  (\ref{eq 6.1-case2})
 and it verifies that
 $$v_0\leq v_{r_1,\infty}\leq v_0+\bar v\quad {\rm in}\ \, \R^N\setminus\{0\},$$
 which implies that
 $$\lim_{|x|\to0^+}v_{r_1,\infty}(x)|x|^{N-2s}(-\ln|x|)^{-m_0}=\cK_s.$$

For $r\leq r_1$,
 $$
\left\{
\begin{array}{lll}
(-\Delta)^s v_{r,n}=   v_{r,n-1}^{p^*}  \quad\,  {\rm in}\ \, B_{r} \setminus\{0\},
\\[2mm]
\qquad \ \ \, v_{r,n}=v_0 \quad\ \ \quad {\rm in}\ \,  \R^N\setminus B_{r},\\[2mm]
\displaystyle \lim_{|x|\to0^+} v_{r,n}(x)|x|^{N-2s}=0
\end{array}\right.
$$
has a unique solution   $v_{r,n}$.  By the comparison principle, there holds
$$v_{r,n}\leq v_{n}\quad {\rm in}\ \,  \R^N\setminus\{0\}$$
and then  the limit of $\{v_{r,n}\}$ exists, denoting $v_{r,\infty}$,
$$v_{r,\infty}\leq v_{r_1,\infty}\quad {\rm in}\ \,  \R^N\setminus\{0\}.$$

 Thanks to the nonnegative property of $v_0$,   the function $v_{r_1,\infty}$,  the solution of  (\ref{eq 6.1-case2}),  is  a super solution of    (\ref{eq 6.1-case1}) for $r\in(0, r_1]$.  Moreover, $v_{r_1,0}-v_{r_1,0}(r_1)$ is a sub solution of    (\ref{eq 6.1-case1}).    Repeat above arguments,
(\ref{eq 6.1-case1}) has a solution $u_{r_1}$ such that 
$$v_{r_1,0} -v_{r_1,0}(r_1)\leq  u_{r_1} \leq v_{r_1,\infty}  \leq v_{r_1,0}  +\bar v\quad{\rm in}\ \,   B_{r_1}\setminus\{0\}.$$
 
 For $0<r\leq r_1$,  comparison principle implies that  
  if $v_{r,1}\leq v_{r_1,1}$ and inductively, we obtain tha t
  $v_{r,\infty}\leq v_{r_1,\infty}$ and  the order could be kept, that is,
 then
  $$u_{r}  \leq u_{r_1} \quad {\rm for\ }\ r\leq r_1.$$
  This completes the proof.\hfill$\Box$\medskip

\begin{corollary}\label{cr existence 2}
 Let  $N>2$,   then for $r\in(0,\frac1{e^2}]$
 problem
  \begin{equation}\label{eq 6.1-case1 s=1}
\left\{
\begin{array}{lll}
 -\Delta  u=   u^{\frac{N}{N-2}}  \quad &{\rm in}\ \, B_{r} \setminus\{0\},
\\[2mm]
\qquad \ \ \, u=0&{\rm on}\ \,  \partial B_{r}
\end{array}\right.
 \end{equation}
 has a positive singular solution $u_r$ verifying for some $m>0$
 $$-\cK_1 r^{2-N}(-\ln r)^{-\frac{N-2}{2}} \leq   u_r(x) - \cK_1 |x|^{2-N}(-\ln|x|)^{\frac{N-2}{2}}\leq  c_{27} |x|^{2-N}(-\ln|x|)^{-\frac{N}{2}},$$
where  $c_{27}>0$ is independent of $r$. 

Moreover, the mapping $r\in(0,\frac1{e^2}] \mapsto u_r$ is increasing.
 \end{corollary}
 {\bf Proof.}  Let
 $$v_0(x)=  |x|^{2-N}(-\ln|x|)^{\frac{2-N}{2}}\quad {\rm in}\ \, B_{\frac1{e^2}}\setminus\{0\}  $$
 and direct computation shows that
  \begin{eqnarray*}
 -\Delta v_0 &= &(-\ln|x|)^{\frac{2-N}{2}}(-\Delta)|x|^{2-N}
 \\[1.5mm]&&-2\nabla|x|^{2-N}\cdot \nabla (-\ln|x|)^{\frac{2-N}{2}}-|x|^{2-N}(-\Delta)(-\ln|x|)^{\frac{2-N}{2}}
 \\[1.5mm]&=&-\frac{(N-2)^2}2|x|^{-N}(-\ln|x|)^{\frac{-N}{2}}-\frac{N(N+2)}{4}|x|^{-N}(-\ln|x|)^{\frac{-N-2}{2}},
\end{eqnarray*}
which gives a sub solution of
\begin{equation}\label{eq 6.1-loc}
\left\{
\begin{array}{lll}
 -\Delta u=   u^{\frac{N}{N-2}}  \quad\,  {\rm in}\ \, B_{r} \setminus\{0\},
\\[2mm]
\quad \ \   u=v_0 \quad\ \ \quad {\rm on}\ \,  \partial B_{r}.
\end{array}\right.
\end{equation}
for $r\in(0, \frac1{e^2})$. 

 By a similar iterative procedure  of  Proposition \ref{existence 2}, we can construct 
 a desired solution $u_r$. \hfill$\Box$\medskip

 \subsection{Infinitely many solutions}
 \begin{lemma}\label{lm 7.1}
 Let $u_r$ be the solution of (\ref{eq 6.1-case1})
 obtained in Proposition \ref{existence 2}, then 
 $$\Big|u_r(x)-\cK_s|x|^{2s-N}(-\ln|x|)^{m_0}  \Big|\leq c_{29}|x|^{2s-N}(-\ln|x|)^{(m_0-1)p^*+1} \quad {\rm for}\ \, x\in B_{\frac{r}2}\setminus\{0\}, $$
 where $c_{29}>0$ depends on $r$ and
 $$(m_0-1)p^*+1<m_0-1.$$
 \end{lemma}
 \noindent{\bf Proof. }
 From Proposition  \ref{existence 2},
  problem  (\ref{eq 6.1-case1})
 has a positive singular solution $u_r$ verifying
 \begin{equation}\label{eq 6.1-0}
\Big| u_r(x)-\cK_s|x|^{2s-N}(-\ln|x|)^{m_0} \Big|\leq   c_{30} |x|^{2s-N}(-\ln|x|)^{m_0-1} \quad {\rm for}\ \, x\in B_{\frac{r}2}\setminus\{0\},
  \end{equation}
  where $c_{30}>0$ depends on $r$.

  Then we have for $x\in B_{\frac{r}2}\setminus\{0\}$:
  $$\Big|u_{r}^{p^*}-\cK_s^{p^*} |x|^{-N}(-\ln|x|)^{m_0-1}\Big|\leq c_{31}|x|^{-N}(-\ln|x|)^{(m_0-1)p^*}  $$
  and
  \begin{eqnarray*}
 \Big| (-\Delta)^s\big( u_r-\cK_s w_{m_0}\big)(x)\Big|&=& \Big|u_r^{p^*}(x)-\cK_s^{p^*} |x|^{-N} (-\ln|x|)^{m_0-1}  \Big|
 \\[1.5mm]&\leq &c_{31}|x|^{-N}(-\ln|x|)^{(m_0-1)p^*},
   \end{eqnarray*}
which, by Lemma \ref{lm 2.2-upper},  implies that
$$\Big|u_r-\cK_s w_{m_0}\Big|\leq c_{32}|x|^{2s-N}(-\ln|x|)^{(m_0-1)p^*+1}.  $$
 This completes the proof.  \hfill$\Box$\medskip

\noindent{\bf Proof of Theorem \ref{teo 2}. }
 Lemma  \ref{lm 7.1}  shows that  for any $r\leq r_1$,
  problem  (\ref{eq 6.1-case1})
 has a positive singular solution $u_{r}$ verifying that for $0<|x|<r $
 \begin{equation}\label{eq 6.1-0}
 \Big|  u_{r}(x)-\cK_s|x|^{2s-N}(-\ln|x|)^{m_0} \Big| \leq  c_{29} |x|^{2s-N}(-\ln|x|)^{(m_0-1)p^*+1}.
  \end{equation}

Next we perform the following scaling: 
$$u_l(x)=l^{2s-N} u_r(l^{-1} x)\quad {\rm for}\ \, x\in\R^N\setminus\{0\}.$$ 
Then $u_l$ is a singular solution of
 \begin{equation}\label{eq 6.1-l}
\left\{
\begin{array}{lll}
(-\Delta)^s u=   u^{p^*}  \quad &{\rm in}\ \, B_{lr} \setminus\{0\},
\\[2mm]
\qquad \ \ \, u=0&{\rm in}\ \,  \R^N\setminus B_{lr}.
\end{array}\right.
 \end{equation}
 Note that for $|x|>0$ sufficiently small,
  \begin{eqnarray*}
  \cK_s|x|^{2s-N}(\ln l-\ln|x|)^{m_0}&=& \cK_s|x|^{2s-N}(- \ln|x|)^{m_0}
 +\cK_s m_0 (\ln l) \,|x|^{2s-N}(- \ln|x|)^{m_0-1}
 \\[1.5mm]&& +O(1)|x|^{2s-N}(- \ln|x|)^{m_0-2}.
  \end{eqnarray*}
  From (\ref{eq 6.1-0}), we have that as $|x|\to0^+$
  \begin{eqnarray*}
  u_l(x)-\cK_s|x|^{2s-N}(-\ln|x|)^{m_0}&=&\cK_s m_0 (\ln l) \,|x|^{2s-N}(- \ln|x|)^{m_0-1}
  \\[1.5mm]&&+O(1)|x|^{2s-N}(- \ln|x|)^{\max\{m_0-2,(m_0-1)p^*+1\}},
   \end{eqnarray*}
where $\max\{m_0-2,(m_0-1)p^*+1\}<m_0-1$.

If we choose $l>0$ such that $rl=1$, then 
\begin{equation}\label{eq 6.1-1}
\left\{
\begin{array}{lll}
(-\Delta)^s u=   u^{p^*}  \quad &{\rm in}\ \, B_1 \setminus\{0\},
\\[2mm]
\qquad \ \ \, u=0&{\rm in}\ \,  \R^N\setminus B_1
\end{array}\right.
 \end{equation}
has a sequence  of solutions $\{u_{l}\}_{l}$ satisfy that
 $$u_{l}(x)-\cK_s|x|^{2s-N}(-\ln|x|)^{m_0}=\cK_s m_0 (\ln l)\,|x|^{2s-N}(- \ln|x|)^{m_0-1}(1+o(1))\quad {\rm for}\ \,  |x|\to 0^+,$$
 where $m_0-1=-\frac{N}{2s}$.
This completes the proof. \hfill$\Box$\bigskip

 \noindent{\bf Proof of Corollary \ref{cr 2}. }
 From Corollary \ref{cr existence 2}, 
\begin{equation}\label{eq 6.1-loc}
\left\{
\begin{array}{lll}
 -\Delta u=   u^{\frac{N}{N-2}}  \quad\,  {\rm in}\ \, B_{r} \setminus\{0\},
\\[2mm]
\quad \ \   u=0 \quad\ \ \quad {\rm on}\ \,  \partial B_{r}.
\end{array}\right.
\end{equation}
 has a solution  $u_r$ for $r\in(0, \frac1{e^2})$ such that 
 $$-\cK_1 r^{2-N}(-\ln r)^{-\frac{N-2}{2}} \leq   u_r(x) - \cK_1 |x|^{2-N}(-\ln|x|)^{\frac{N-2}{2}}\leq  c_{27} |x|^{2-N}(-\ln|x|)^{-\frac{N}{2}},$$
where  $c_{27}>0$ is independent of $r$. 

 For $r\in(0,\frac1{e^2})$, we do the scaling
 $$u_l(x)=l^{N-2}u_r(l^{-1}x),\quad x\in B_{rl}\setminus\{0\}.$$
 The rest of the proof  is omitted here.   \hfill$\Box$

\bigskip
\bigskip
 % \newpage
 \appendix
 \section{Appendix: The constant $\cC'_s(0)$ and $\cK_s$}
We show that
\begin{equation}\label{sub 6.1-app}
  \cC'_s(0)= -2^{2\s-1}\frac{\Gamma(\frac{N}{2})\Gamma(\s)}{ \Gamma(\frac{N-2\s}{2})}.
 \end{equation}
 In fact,
 \begin{eqnarray*}
\ln \cC_s(\tau)&= & 2s \ln 2 + \ln \Gamma(\frac{N+\tau}{2}) +\ln \Gamma(\frac{2\s-\tau}{2}) -\ln \Gamma(-\frac{\tau}{2})-\ln \Gamma(\frac{N-2\s+\tau}{2})
 \end{eqnarray*}
 and then
  \begin{eqnarray*}
  \cC_s'(\tau)&= &\cC_s(\tau) \Big(  \frac12 \psi(\frac{N+\tau}{2}) -\frac12 \psi(\frac{2\s-\tau}{2}) +\frac12 \psi  (-\frac{\tau}{2})-\frac12 \psi(\frac{N-2\s+\tau}{2}) \Big)
  \\[1.5mm]&=&\frac{2^{2\s-1}}{\Gamma(-\frac{\tau}{2})} \frac{\Gamma(\frac{N+\tau}{2})\Gamma(\frac{2\s-\tau}{2})}{\Gamma(\frac{N-2\s+\tau}{2})}\left(  \psi(\frac{N+\tau}{2}) -\psi(\frac{2\s-\tau}{2}) -\psi(\frac{N-2\s+\tau}{2})\right)
   \\[1.5mm]&&  + 2^{2\s-1} \frac{\Gamma(\frac{N+\tau}{2})\Gamma(\frac{2\s-\tau}{2})}{\Gamma(\frac{N-2\s+\tau}{2})} \frac{\psi  (-\frac{\tau}{2})}{\Gamma(-\frac{\tau}{2})}  ,
 \end{eqnarray*}
 where  $\psi $ is the Digamma function, i.e. $\psi(t)=\frac{\Gamma'(t)}{\Gamma(t)}$.
Note that the term $$  \psi(\frac{N+\tau}{2}) -\psi(\frac{2\s-\tau}{2}) -\psi(\frac{N-2\s+\tau}{2})$$
 is uniformly bounded as $s\to1^+$

By properties of Gamma function and Digamma function,  we have that
\begin{eqnarray*}
\frac1{\Gamma(-\frac{\tau}{2})}=-\frac{\tau}{2}\frac{1 }{\Gamma(1-\frac{\tau}{2}) }
 \end{eqnarray*}
 and
 \begin{eqnarray*}
\psi  (-\frac{\tau}{2}) =\psi  (1-\frac{\tau}{2})  +\frac{\tau}{2}
 \end{eqnarray*}
then
 $$\lim_{\tau\to0}    \frac{\psi  (-\frac{\tau}{2})}{\Gamma(-\frac{\tau}{2})}=-1 $$
 and
\begin{eqnarray*}
\lim_{\tau\to0}  \cC_s'(\tau)&= & 2^{2\s-1}  \lim_{\tau\to0}\frac{\Gamma(\frac{N+\tau}{2})\Gamma(\frac{2\s-\tau}{2})}{\Gamma(\frac{N-2\s+\tau}{2})} \lim_{\tau\to0}  \frac1{\Gamma(-\frac{\tau}{2})}\Big(\psi(\frac{N+\tau}{2}) -\psi(\frac{2\s-\tau}{2})
\\[2mm]&&-\psi(\frac{N-2\s+\tau}{2})\Big) + 2^{2\s-1} \lim_{\tau\to0}\frac{\Gamma(\frac{N+\tau}{2})\Gamma(\frac{2\s-\tau}{2})}{\Gamma(\frac{N-2\s+\tau}{2})} \lim_{\tau\to0}    \frac{\psi  (-\frac{\tau}{2})}{\Gamma(-\frac{\tau}{2})}
\\[2mm]&=& -2^{2\s-1} \frac{\Gamma(\frac{N}{2})\Gamma(\frac{2\s }{2})}{\Gamma(\frac{N-2\s }{2})}.
 \end{eqnarray*}

Next we prove that
\begin{equation}\label{sub 6.1-app1}
 \lim_{s\to1^-}\cK_s=\Big(\frac{(N-2)^2}{2}\Big)^{\frac{N-2}{2}}.
 \end{equation}

Notice that
\begin{eqnarray*}
\cK_s&= &  \Big(2^{2\s-1} \frac{\Gamma(\frac{N}{2})\Gamma(\frac{2\s }{2})}{\Gamma(\frac{N-2\s }{2})}\frac{N-2s}{2s}\Big)^{\frac{N-2s}{2s}},
 \end{eqnarray*}
 where
\begin{eqnarray*}
\lim_{s\to1^-} \Big(2^{2\s-1} \frac{\Gamma(\frac{N}{2})\Gamma(\frac{2\s }{2})}{\Gamma(\frac{N-2\s }{2})}
 \frac{N-2s}{2s} \Big)&=&(N-2)\frac{\Gamma(\frac{N}{2})}{\Gamma(\frac{N-2}{2}) }
 \\[2mm]&=&(N-2)\frac{\frac{N-2}{2} \Gamma(\frac{N-2}{2})}{\Gamma(\frac{N-2}{2}) }=\frac{(N-2)^2}{2}.
 \end{eqnarray*}

\section{Appendix: Radial symmetry}

 \setcounter{equation}{0}

 Our method of moving planes is motivated by \cite{FW} to deal with the solution $u$ has possibly singular at the origin.  For the moving planes of integral equations, we refer to \cite{CLO1,CLO}.
For singular solutions, we use a direct moving plane method from  \cite{FW} and we need the following variant Maximum Principle for small domain.
\begin{lemma}\label{pr 2-a}\cite[Corollary 2.1]{FW} Let $O$  be an
open and bounded subset of $\R^N$. Suppose that
$\varphi:\Omega\to\R$ is in $L^\infty(O)$ and
$w\in L^\infty(\R^N)$ is a classical solution of
\begin{equation}\label{eq a1-a}
\left\{ \arraycolsep=1pt
\begin{array}{lll}
-(-\Delta)^s  w(x)\le \varphi(x)w(x),\ \ \ \ &
x\in O,\\[2mm]
\qquad\quad\ w(x)\ge 0,& x\in  O^c.
\end{array}
\right.
\end{equation}
Then there is $\delta>0$ such that whenever $|O^-|\le \delta$,
$w$ has to be  non-negative in $O$, where $O^-=\{x\in O\ | \ w(x)<0\}$.
 \end{lemma}

Now, we use the moving plane method to show the radial symmetry and monotonicity of positive solutions to equation (\ref{eq 5.1}).
  For simplicity, we denote
\begin{equation}\label{d1}\Sigma_\lambda=\{x=(x_1,x')\in \mathcal{B}_1\  |\  x_1>\lambda\},\end{equation}
\begin{equation}\label{d3}u_\lambda(x)=u(x_\lambda) \quad\mbox{and}\quad  w_\lambda(x)=u_\lambda(x)-u(x),\end{equation}
where $\lambda\in (0,1)$ and $x_\lambda=(2\lambda-x_1,x')$  for
$x=(x_1,x')\in\R^N$ and  $\mathcal{B}_1=B_1 \setminus\{0\}$.  For any subset $A$ of $\R^N$, we write
$A_\lambda=\{x_\lambda:\, x\in A\}$.

 On the contrary, suppose that $\Sigma_\lambda^-=\{x\in \Sigma_\lambda\ |\ w_\lambda(x)<0\}\not=\emptyset$ for $\lambda\in(0,1)$.
Let us define
\begin{equation}\label{eq c}
w_\lambda^+(x)=\left\{ \arraycolsep=1pt
\begin{array}{lll}
w_\lambda(x),\ \ \ \ &
x\in \Sigma_\lambda^-,\\[1mm]
0,& x\in \R^N\setminus \Sigma_\lambda^-
\end{array}
\right.
\end{equation}
and
\begin{equation}\label{eq 4.01}
w_\lambda^-(x)=\left\{ \arraycolsep=1pt
\begin{array}{lll}
0,\ \ \ \ &
x\in \Sigma_\lambda^-,\\[1mm]
w_\lambda(x),\ \ \ \ & x\in \R^N\setminus \Sigma_\lambda^-.
\end{array}
\right.
\end{equation}
Hence, $w_\lambda^+(x)=w_\lambda(x)-w_\lambda^-(x)$ for all $x\in\R^N.$ It is obvious that
 $(2\lambda,0,\cdots,0)\not\in \Sigma_\lambda^-$, since $\displaystyle \lim_{|x|\to0^+}u(x)=+\infty$.

\begin{lemma}\label{lm 13-08-1}
Assume that $\Sigma_\lambda^-\not=\emptyset$  for  $ 0<\lambda<1$, then
\begin{equation}\label{claim1}
 (-\Delta)^s w_\lambda^-(x)\le0,\ \ \ \ \forall\ x\in \Sigma_\lambda^-.
 \end{equation}

\end{lemma}
{\bf Proof.}  By direct computation, for
 $x\in \Sigma_\lambda^-$, we have
 \begin{eqnarray*}
 (-\Delta)^s w_\lambda^-(x)
&=&\int_{\R^N}\frac{w_\lambda^-(x)-w_\lambda^-(z)}{|x-z|^{N+2s}}dz
=-\int_{\R^N\setminus\Sigma_\lambda^-}\frac{w_\lambda(z)}{|x-z|^{N+2s}}dz
\\&=&-\int_{(\mathcal{B}_1\setminus(\mathcal{B}_1)_\lambda) \cup ((\mathcal{B}_1)_\lambda\setminus \mathcal{B}_1)}\frac{w_\lambda(z)}{|x-z|^{N+2s}}dz
\\&&-\int_{(\Sigma_\lambda\setminus\Sigma_\lambda^-) \cup (\Sigma_\lambda\setminus\Sigma_\lambda^-)_\lambda}\frac{w_\lambda(z)}{|x-z|^{N+2s}}dz
%\\&&
-\int_{(\Sigma_\lambda^-)_\lambda}\frac{w_\lambda(z)}{|x-z|^{N+2s}}dz
\\&=&-I_1-I_2-I_3.
\end{eqnarray*}
We estimate these integrals separately. Since $u=0\ {\rm in} \ (\mathcal{B}_1)_\lambda\setminus \mathcal{B}_1$ and
$u_\lambda=0\ {\rm in} \ \mathcal{B}_1\setminus (\mathcal{B}_1)_\lambda$,
then
\begin{eqnarray*}
I_1
&=&\int_{(\mathcal{B}_1\setminus(\mathcal{B}_1)_\lambda) \cup ((\mathcal{B}_1)_\lambda\setminus \mathcal{B}_1)}\frac{w_\lambda(z)}{|x-z|^{N+2s}}dz
\\&=&\int_{(\mathcal{B}_1)_\lambda\setminus \mathcal{B}_1}\frac{u_\lambda(z)}{|x-z|^{N+2s}}dz-\int_{ \mathcal{B}_1\setminus(\mathcal{B}_1)_\lambda}\frac{u(z)}{|x-z|^{N+2s}}dz
\\&=&\int_{(\mathcal{B}_1)_\lambda\setminus \mathcal{B}_1}u_\lambda(z)\Big(\frac{1}{|x-z|^{N+2s}}-\frac{1}{|x-z_\lambda|^{N+2s}}\Big)dz\ge 0,
\end{eqnarray*}
since $u_\lambda\geq0$ and $|x-z_\lambda|>|x-z|$ for all  $x\in \Sigma_\lambda^-$ and $z\in (\mathcal{B}_1)_\lambda\setminus \mathcal{B}_1.$
%we have  $|x-z_\lambda|>|x-z|$.
%Therefore, by $u_\lambda\geq0\ in \ (B_1)_\lambda\setminus B_1,$ we have
%$$I_1\geq0.$$

In order to decide the sign of $I_2$ we observe that $w_\lambda(z_\lambda)=-w_\lambda(z)$ for any $z\in\R^N$. Then,
\begin{eqnarray*}
I_2
&=&\int_{(\Sigma_\lambda\setminus\Sigma_\lambda^-) \cup(\Sigma_\lambda\setminus\Sigma_\lambda^-)_\lambda}\frac{w_\lambda(z)}{|x-z|^{N+2s}}dz
\\&=&\int_{\Sigma_\lambda\setminus\Sigma_\lambda^-}\frac{w_\lambda(z)}{|x-z|^{N+2s}}dz
+\int_{\Sigma_\lambda\setminus\Sigma_\lambda^-}\frac{w_\lambda(z_\lambda)}{|x-z_\lambda|^{N+2s}}dz
\\&=&\int_{\Sigma_\lambda\setminus\Sigma_\lambda^-}w_\lambda(z)(\frac{1}{|x-z|^{N+2s}}-\frac{1}{|x-z_\lambda|^{N+2s}})dz
\\&\ge& 0,
\end{eqnarray*}
since  $w_\lambda\ge0$ in $\Sigma_\lambda\setminus\Sigma_\lambda^-$ and $|x-z_\lambda|>|x-z|$ for all $x\in \Sigma_\lambda^-$ and $z\in
\Sigma_\lambda\setminus\Sigma_\lambda^-.$
%we have
%$$I_2\geq0.$$

Finally, since
  $w_\lambda(z)<0$ for $z\in \Sigma_\lambda^-$, we deduce
\begin{eqnarray*}
I_3
&=&\int_{(\Sigma_\lambda^-)_\lambda}\frac{w_\lambda(z)}{|x-z|^{N+2s}}dz
=\int_{\Sigma_\lambda^-}\frac{w_\lambda(z_\lambda)}{|x-z_\lambda|^{N+2s}}dz
\\&=&-\int_{\Sigma_\lambda^-}\frac{w_\lambda(z)}{|x-z_\lambda|^{N+2s}}dz
\ge0.
\end{eqnarray*}
The proof is complete.\hfill$\Box$

\medskip

Now we are ready to prove Theorem \ref{teo 5.1-1}.\medskip

\noindent{\bf Proof of  Theorem \ref{teo 5.1-1}.}
Our purpose is to show the radial symmetry and decreasing monotonicity and we divide the proof into four steps.\smallskip

\noindent \emph{Step 1:} We prove that if $\lambda$ is close to $1$,  then  $w_\lambda>0$  in $\Sigma_\lambda$.
First we show that  $w_\lambda\geq0$  in $\Sigma_\lambda$, i.e. $\Sigma^-_\lambda$ is empty.
By contradiction, we assume that
$\Sigma^-_\lambda\not=\emptyset$.
Now we apply  \equ{claim1} and linearity of the fractional Laplacian to   obtain that, for $ x\in\Sigma_\lambda^-,$
\begin{equation}\label{eq e}
(-\Delta)^s w_\lambda^+(x)\ge (-\Delta)^s w_\lambda(x)
=(-\Delta)^s u_\lambda(x)-(-\Delta)^s u(x).
\end{equation}
Combining    with (\ref{eq e}) and  (\ref{eq c}), for $x\in\Sigma_\lambda^-$, we have
\begin{eqnarray*}
(-\Delta)^s w_\lambda^+(x) &\geq&(-\Delta)^s
u_\lambda(x)-(-\Delta)^s u(x)
\\&=&-u_\lambda^{p^*}(x) + u^{p^*}(x)
 = - \varphi(x)w_\lambda^+(x),
\end{eqnarray*}
where
$$
   \varphi(x)= \frac{ (u_\lambda(x))^{p^*}-(u(x))^{p^*}}{u_\lambda(x)-u(x)},\quad\forall\, x\in\Sigma_\lambda^-.
$$

For $x\in \Sigma_\lambda^-\subset \Sigma_\lambda\subset \R^N\setminus {B_\lambda(0)}$, $u_\lambda(x)<u(x)$. Moreover,
there exists $M_\lambda>0$ such that
$$\norm{u}_{L^\infty(\R^N\setminus {B_\lambda(0)})}\le M_\lambda.$$
Due to $h\in C^1(\R_+)$,   there exists $c_{31}>0$ dependent of $\lambda$ such that
\begin{equation}\label{app 1}
\norm{\varphi}_{ L^\infty(\Sigma_\lambda^-)}\le c_{31}.
\end{equation}
Note that  $M_\lambda\to \infty$ as $\lambda\to 0$, since $\displaystyle \lim_{|x|\to0^+}u(x)=\infty$.\smallskip

Therefore, for  $x\in\Sigma_\lambda^-$ and then
$$%\begin{equation}\label{eq f}
-(-\Delta)^s  w_\lambda^+(x)\le \varphi(x)w_\lambda^+(x), \quad
\forall\, x\in\Sigma_\lambda^-.
$$%\end{equation}
Moreover, $w_\lambda^+=0$ in $(\Sigma_\lambda^-)^c$. Choosing $\lambda\in (0,1)$ close enough to $1$ we have $|\Sigma_\lambda^-|$ is small
and we apply Lemma \ref{pr 2-a} to obtain that
$$w_\lambda=w_\lambda^+\geq0\ \ \ \ \mbox{in} \ \ \Sigma_\lambda^-,$$
which is impossible. Thus,
$$w_\lambda\geq0\ \ \ \mbox{in}\ \ \Sigma_\lambda.$$

Now we claim that for $0<\lambda<1$, if
 $w_\lambda\ge0$ and $w_\lambda\not\equiv0$
 in $\Sigma_\lambda$, then  $w_\lambda>0$
 in $\Sigma_\lambda$. Assuming the claim is true, we complete the proof. Since the function $u$ is positive in $B_1 $ and
 $u=0$ on $\partial B_1 $, $w_\lambda$ is positive on $\partial B_1 \cap \partial \Sigma_\lambda$ and then  $w_\lambda\not\equiv0$ in $\Sigma_\lambda$.

Now we prove the claim.  Suppose on the contrary that there exists $x_0\in \Sigma_\lambda$ such that
$w_\lambda(x_0)=0,$ i.e.  $u_\lambda(x_0)=u(x_0)$. Then
\begin{eqnarray}\label{eq 7}
(-\Delta)^s w_\lambda(x_0)= (-\Delta)^s
u_\lambda(x_0)-(-\Delta)^s u(x_0)
=0.
\end{eqnarray}

On the other hand, let
$K_\lambda=\{(x_1,x')\in\R^N\ |\  x_1>\lambda\}$. Noting
$w_\lambda(z_\lambda)=-w_\lambda(z)$ for any $z\in\R^N$ and $w_\lambda(x_0)=0$, we deduce
\begin{eqnarray*}
(-\Delta)^s w_\lambda(x_0)
&=&-\int_{K_\lambda}\frac{w_\lambda(z)}{|x_0-z|^{N+2s}}dz-\int_{\R^N\setminus K_\lambda}\frac{w_\lambda(z)}{|x_0-z|^{N+2s}}dz
\\&=&-\int_{K_\lambda}\frac{w_\lambda(z)}{|x_0-z|^{N+2s}}dz-\int_{K_\lambda}\frac{w_\lambda(z_\lambda)}{|x_0-z_\lambda|^{N+2s}}dz
\\&=&-\int_{K_\lambda}w_\lambda(z)\Big(\frac{1}{|x_0-z|^{N+2s}}-\frac{1}{|x_0-z_\lambda|^{N+2s}}\Big)dz.
\end{eqnarray*}
The fact $|x_0-z_\lambda|>|x_0-z|$ for $z\in K_\lambda$ ,
$w_\lambda(z)\ge0$ and $w_\lambda(z)\not\equiv0$ in $K_\lambda$ yield
$$%\begin{equation}\label{ff}
(-\Delta)^s w_\lambda(x_0)<0,$$%\end{equation}
which contradicts  (\ref{eq 7}), completing the proof of the claim.

\medskip

\noindent \emph{Step 2:} We prove $\lambda_0:=\inf\{\lambda\in(0,1)\ |\  w_\lambda>0\ \ \rm{in}\ \ \Sigma_\lambda\}=0$. Were it not true, we would have $\lambda_0>0$. Hence,
 $w_{\lambda_0}\ge0$ in $\Sigma_{\lambda_0}$ and
$w_{\lambda_0}\not\equiv0$ in $\Sigma_{\lambda_0}$. The claim in Step 1 implies
 $w_{\lambda_0}>0$ in $\Sigma_{\lambda_0}$.

Next we claim that if $w_\lambda>0$ in $\Sigma_\lambda$ for
$\lambda\in(0,1)$, then there exists $\epsilon\in(0,\lambda/4)$ such
that $w_{\lambda_\epsilon}>0$ in $\Sigma_{\lambda_\epsilon}$, where
$\lambda_\epsilon=\lambda-\epsilon>3\lambda/4$. This claim directly implies that $\lambda_0=0$, which contradicts to the fact $\lambda_0>0$.

Now we  prove the claim.
Let $D_\mu=\{x\in\Sigma_\lambda\ | \ dist(x,\partial\Sigma_\lambda)\ge \mu\}$ for $\mu>0$ small. Since
$w_\lambda>0$ in $\Sigma_\lambda$ and $D_\mu$ is compact, there
exists $\mu_0>0$ such that $w_\lambda\ge \mu_0$ in $D_\mu$. By
the continuity of $w_\lambda(x)$, for $\epsilon>0$ small
enough and  $\lambda_\epsilon=\lambda-\epsilon,$ we have that
$w_{\lambda_\epsilon}(x)\ge0\ \ \rm{in}\ \ D_\mu$. Therefore,
$\Sigma_{\lambda_\epsilon}^-\subset
\Sigma_{\lambda_\epsilon}\setminus D_\mu$ and
$|\Sigma_{\lambda_\epsilon}^-|$ is small if $\epsilon$ and $\mu$ are small.
Using \equ{claim1} and proceeding as in  Step 1, we have for all $x\in \Sigma_{\lambda_\epsilon}^-$
that \begin{eqnarray*}
(-\Delta)^s w_{\lambda_\epsilon}^+(x) &=&(-\Delta)^s
u_{\lambda_\epsilon}(x)-(-\Delta)^s u(x)-(-\Delta)^s
w_{\lambda_\epsilon}^-(x)
\\&\ge& (-\Delta)^s u_{\lambda_\epsilon}(x)-(-\Delta)^s u(x)
\\&=& \varphi(x)w_{\lambda_\epsilon}^+(x).
\end{eqnarray*}
By (\ref{app 1}), if $\lambda_\epsilon>3\lambda/4$,
$\varphi(x)$ is controlled by some constant dependent of $\lambda$.

Since $w_{\lambda_\epsilon}^+=0$ in
$(\Sigma_{\lambda_\epsilon}^-)^c$ and
$|\Sigma_{\lambda_\epsilon}^-|$ is small, for  $\epsilon$ and $\mu$
small,  Lemma \ref{pr 2-a}  implies that  $w_{\lambda_\epsilon}\ge0$ in
$\Sigma_{\lambda_\epsilon}$.  Combining with $\lambda_\epsilon>0$ and
$w_{\lambda_\epsilon}\not\equiv0$ in $\Sigma_{\lambda_\epsilon}$, we obtain
 $w_{\lambda_\epsilon}>0$
 in $\Sigma_{\lambda_\epsilon}$.  The proof of the claim is finished.

\medskip

\noindent \emph{Step 3:} By Step 2, we have $\lambda_0=0$, which
implies that $u(-x_1,x')\ge u(x_1,x')$ for $x_1\ge0.$
Using the same argument from the other side, we conclude that $u(-x_1,x')\le u(x_1,x')$ for $ x_1\ge0$ and then
$u(-x_1,x')= u(x_1,x')$ for $x_1\ge0.$ Repeating this procedure in all directions we see that $u$ is radially symmetric.

\medskip

Finally, we prove $u(r)$ is strictly decreasing in  $r\in (0,1)$. Let us consider  $0<x_1<\widetilde{x}_1<1$
and let  $\lambda=\frac{x_1+\widetilde{x}_1}{2}$. As proved above   we have
$$w_\lambda(x)>0\ \ \mbox{for}\ \ x\in\Sigma_{\lambda}.$$
Then
\begin{eqnarray*}
0<w_\lambda(\widetilde{x}_1,0,\cdots,0)
&=&u_\lambda(\widetilde{x}_1,0,\cdots,0)-u(\widetilde{x}_1,0,\cdots,0)
\\&=&u(x_1,0,\cdots,0)-u(\widetilde{x}_1,0,\cdots,0),
\end{eqnarray*}
i.e
$u(x_1,0,\cdots,0)>u(\widetilde{x}_1,0,\cdots,0).$
From the radial symmetry of $u$ and decreasing in the direction $\frac{x}{|x|}$,  we can conclude the monotonicity of $u$.
 \hfill$\Box$\medskip

\setcounter{equation}{0}
  \section{Non-existence in the whole space }

Let $\cL_\mu^s$ be the fractional Hardy operator  defined by
$$ \cL_\mu^s  = (-\Delta)^s    +\frac{\mu}{|x|^{2s}} $$
for $\s\in(0,1)$   and $$\mu\geq  \mu_0 :=-2^{2\s}\frac{\Gamma^2(\frac{N+2\s}4)}{\Gamma^2(\frac{N-2\s}{4})}.$$

It is shown in \cite{CW} that  for $\mu\geq \mu_0$ the equation
$$\mathcal{ L}_\mu^s u=0\quad{\rm in}\ \ \R^N\setminus \{0\}$$
 has two distinct radial solutions
 $$\Phi_{s,\mu}(x)=\left\{\arraycolsep=1pt
\begin{array}{lll}
 |x|^{\tau_-(s,\mu)}\quad
   &{\rm if}\quad \mu>\mu_0\\[1.5mm]
 \phantom{   }
|x|^{-\frac{N-2s}{2}}\ln\left(\frac{1}{|x|}\right) \quad  &{\rm   if}\quad \mu=\mu_0
 \end{array}
 \right.\qquad \ {\rm and}\quad\ \ \Gamma_{s,\mu}(x)=|x|^{\tau_+(s,\mu)},$$
where
$\tau_-(s,\mu)  \leq  \tau_+(s,\mu)$ verifies that
 \begin{eqnarray*}
 &\tau_-(\s,\mu)+\tau_+(\s,\mu) =2\s-N \quad\  {\rm for\ all}\ \   \mu \ge \mu_0,\\[1.5mm]
&\tau_-(\s,\mu_0)=\tau_+(\s,\mu_0)=\frac{2\s-N}2,\quad\
\tau_-(\s,0)=2\s-N, \quad\ \tau_+(\s,0)=0,  \\[1.5mm]
&\displaystyle\lim_{\mu\to+\infty} \tau_-(\s,\mu)=-N\quad {\rm and}\quad \lim_{\mu\to+\infty} \tau_+(\s,\mu)=2\s.
 \end{eqnarray*}
For simplicity,  we put $\tau_+=\tau_+(s,\mu)$, $\tau_-=\tau_-(s,\mu)$.

 \begin{theorem}\label{teo 3-l}
Let   $\mu>\mu_0$ and $p\in(1,p^*_\mu]$,
where
$$p^*_{\mu}=  1+\frac{2s }{-\tau_-(s,\mu)}.$$

Then problem
   \begin{equation}\label{eq 1.1-ext l}
(-\Delta)^s u=   u^{p}  \quad  {\rm in}\ \, \R^N\setminus\{0\}
 \end{equation}
    has no positive solution.
\end{theorem}
Note that    the mapping $\mu\in[\mu_0,+\infty)\mapsto p^*_\mu$ is strictly decreasing. Particularly, $p^*_\mu=\frac{N}{N-2s}$ for $\mu=0$ and
$p^*_\mu=\frac{N+2s}{N-2s}$ for $\mu=\mu_0$.
\subsection{Basic estimates}
 \begin{lemma}\label{re 2.1} Let $\mu\ge\mu_0$,  $O$ is a bounded domain containing the origin and nonnegative function  $f\in C^\beta_{loc}(\R^N \setminus \bar O)$ for some $\beta\in(0,1)$.
The homogeneous problem
\begin{equation}\label{eq 1.1 EH}
\left\{\arraycolsep=1pt
\begin{array}{lll}
 \cL^\s_\mu u  \geq  f   \quad
   &{\rm in}\ \ \R^N\setminus\bar O, \\[2mm]
 \phantom{ \cL_\mu^s   }
u\geq 0 \quad &{\rm in}\ \ \,  \bar O
 \end{array}
 \right.
\end{equation}
 has no positive solution if
 $$\lim_{r\to+\infty}\int_{B_r(0)\setminus B_{R_0}(0)}f(x)|x|^{ 2s-\tau_+-N}dx=+\infty.$$
Particularly, the above assumption could be replaced by
$$\liminf_{|x|\to+\infty} f(x)|x|^{2s -\tau_+} >0.$$

\end{lemma}
   \noindent{\bf Proof. }
 Let $u_f$ be positive solution of  (\ref{eq 1.1 EH}) and  denote by
 \begin{eqnarray}\label{kt-1-1}
O^\sharp=\left\{x\in\R^N: \ \frac{x}{|x|^2}\in O\right\}\quad{\rm and}\quad  u^\sharp(x)=|x|^{2s-N}u_f(\frac{x}{|x|^2})\quad{\rm for}\  x\in O^\sharp.
\end{eqnarray}
  Clearly, $(\R^{N}\setminus\{0\})^\sharp=\R^{N}\setminus\{0\}$  and the Kelvin transformation is the following
 \begin{eqnarray}\label{kt-1}
 (-\Delta)^s u^\sharp(x)    =  |x|^{-2s-N} \big((-\Delta)^s u_f\big)\left(\frac{x}{|x|^2}\right)\quad{\rm for}\  x\in O^\sharp.
 \end{eqnarray}
Therefore,  $u^\sharp$ is a super solution of
\begin{equation}\label{eq 1.1 KT}
\cL_{\mu}^s u^\sharp(x) \geq   |x|^{-2s-N} \tilde f (x)  \quad    {\rm in}\ \  O^\sharp,
\end{equation}
where $\tilde  f (x)=f(\frac{x}{|x|^2})$.

 Note that  there exists $r_1>0$ such that
$$B_{r_1}(0)\subset O^\sharp\cup\{0\}.$$
Note that for $r>\frac1{r_1}>0$
 \begin{eqnarray*}
\int_{B_{r_1}(0)\setminus B_\frac1r(0) }f(\frac{x}{|x|^2})|x|^{ -2s-N+\tau_+}dx&=& \int_{B_r(0)\setminus B_{\frac1{r_1}}(0)}f(y)|y|^{ 2s-\tau_+-N}dy
%\\&=&\int_{B_r(0)\setminus B_{\frac1{r_1}}(0)}f(x)|x|^{ N+\tau_-}dx
\\&\to&+\infty\quad{\rm as}\ \,  r\to+\infty.
 \end{eqnarray*}
The contradiction follows by   \cite[Theorem 1.3]{CW}. We complete the proof.
\hfill$\Box$

 \begin{lemma}\label{lm 2.1-ex}
Let  $\mu\ge \mu_0$,  $\theta\in\R$ and $u_0$ be a positive solution of (\ref{eq 1.1-ext}) in $\R^N\setminus\Omega$,
 then
 $$\liminf_{|x|\to+\infty} u_0(x) |x|^{N-2s}>0.$$
 \end{lemma}
 {\bf Proof. } Let $f=Qu^p$, then there exists $\varepsilon_0>0$ and $r_1>r_2>0$  such that
  $$f\geq \varepsilon_0\quad {\rm in}\ \ B_{r_1}\setminus B_{r_2}.$$
 Our problem reduces to
  \begin{equation}\label{eq 2.1-green-ex}
\left\{\arraycolsep=1pt
\begin{array}{lll}
(-\Delta)^s u\geq  f  \quad\       &{\rm in}\ \,   \R^N\setminus \bar \Omega ,
\\[2mm]
\qquad \ \  u\geq 0\quad  \   & {\rm in}\ \     \bar\Omega.%\\[2mm]
%\displaystyle\liminf_{|x|\to+\infty}v(x)|x|^{-\tau_+}=0
  \end{array}
 \right.
\end{equation}

Let
$$
  v^\sharp(x)=|x|^{2s-N}u(\frac{x}{|x|^2})\quad{\rm for}\  x\in O^\sharp,
 $$
and direct computation shows that
$$(-\Delta)^s v^\sharp(x)=|x|^{-2s-N}  f(\frac{x}{|x|^2})\quad{\rm in}\ \ B_1\setminus\{0\}.$$
 Then Maximum principle shows that
$$v^\sharp(x)\geq c_{32} \quad{\rm in}\ \ B_{\frac12} \setminus\{0\},$$
which implies that
$$u_0(\frac{x}{|x|^2})\geq c_{32}|x|^{N-2s}\quad{\rm in}\ \ B_{\frac12}\setminus\{0\}$$
and then
$$u_0(x)\geq c_{32}|x|^{ 2s-N}\quad{\rm in}\ \ \R^N\setminus B_2, $$
where $\tau_++\tau_-=2s-N$.
  \hfill$\Box$\medskip

  \begin{lemma}\label{lm 2.1-ex} Let $2s+\theta>0$, $\tau_0<0$
  \begin{equation}\label{2.00-ex}
 p\in\left[1,\ 1+\frac{2s+\theta}{-\tau_0}\right)
\end{equation}
     and   $\{\tau_j\}_j$ be the sequence generated by
  $$ \tau_j=2s+\theta+p\tau_{j-1}\quad{\rm for}\quad  j=1,2,3\cdots.$$

Then $\{\tau_j\}_j$ is an increasing sequence of numbers and for any $\bar \tau>\tau_0$
there exists $j_0\in\N $  such that
$$
 \tau_{j_0}\ge \bar \tau\quad {\rm and}\quad \tau_{j_0-1}<\bar \tau.
$$
\end{lemma}
 \noindent{\bf Proof. }
For $p\in(0,1+\frac{2s+\theta}{-\tau_0})$, we have that
$$\tau_1-\tau_0=2\alpha+\tau_0(p-1)>0$$
and
$$
\tau_j-\tau_{j-1} = p(\tau_{j-1}-\tau_{j-2})=p^{j-1} (\tau_1-\tau_0),
$$
which imply that the sequence $\{\tau_j\}_j$ is  increasing and our conclusions are obvious.
\hfill$\Box$\medskip

     \begin{lemma}\label{teo 2.2} Assume that $\mu>\mu_0$,
 $O_r=\R^N\setminus\bar B_\frac1r$, the nonnegative function $g\in C^\beta_{loc}(O_r)$ for some $\beta\in(0,1)$ and there exist  $\tau\in(\tau_- ,\tau_+ )$, $c_{33}>0$ and $r_3>0$ such that
 $$ g(x)\geq c_{33}|x|^{\tau-2s}\quad{\rm in}\ \, O_{r_3}.$$

Let  $u_g$ be a positive solution of problem
 \begin{equation}\label{eq 2.1-hom}
\cL_\mu^s u\geq g \quad {\rm in}\ \ O_r,\qquad  u\geq0 \quad  {\rm in}\ \  \R^N\setminus   O_r,
 \end{equation}
then there exists $c_{34}>0$  such that
$$u_g(x)\geq c_{34}|x|^\tau\quad{\rm in}\ \  O_{r_3}.$$

\end{lemma}
{\bf Proof. }  For $\tau\in(\tau_- ,\tau_+ )$, we have that
$$\cL^s_{\mu} |x|^{\tau}=b_s(\tau) |x|^{\tau-2s}\quad  {\rm in}\ \  \R^N\setminus\{0\},$$
where $b_s(\tau)>0$.

In the case $O=B_1\setminus \{0\}$, we use the function
$$w(x)=|x|^{\tau}-|x|^{\tau_+}\quad  {\rm in}\ \  \R^N\setminus\{0\}$$
as a sub solution
$$\cL^s_{\mu} u=c_s(\tau)|x|^{\tau-2s}\quad{\rm in}\ B_1\setminus\{0\},\qquad
u\leq 0\ \ {\rm in} \ \ \R^N\setminus B_1,$$
where $c_s(\tau)>0$.
 Then our argument follows by comparison principle.

 When $O=\R^N\setminus B_1$,  we use the function
$$w(x)=|x|^{\tau}-|x|^{\tau_-}\quad  {\rm in}\ \  \R^N\setminus\{0\}$$
as a sub solution
$$\cL^s_{\mu} u=c_s(\tau)|x|^{\tau-2s}\quad{\rm in}\ \R^N\setminus \bar B_1,\qquad
u\leq 0\ \ {\rm in} \ \  \bar B_1$$
 and the left is standard.\hfill$\Box$

 \subsection{Nonexistence}

 \noindent {\bf Proof of Theorem \ref{teo 3-l}. } By contradiction,
we assume that  (\ref{eq 1.1-ext l}) has a positive solution $u_0$. 
From Lemma \ref{lm 2.1-ex}, we have that
$$u_0(x)\geq d_0 |x|^{2s-N}\quad {\rm in}\ \, \R^N\setminus     B_2.$$

  {\bf Step 1: } We first show the nonexistence of
\begin{equation}\label{eq 1.1-ext-1}
(-\Delta)^s u+\frac{\mu}{|x|^{2s}}=   u^{p}  \quad  {\rm in}\ \, \R^N\setminus\{0\},
 \end{equation}
under the assumption that  $\mu\in[0,\mu_0)$   and $p\in(0, p^*_{\mu})$,
where
\begin{equation}\label{eq 1.1-cr1}
p^*_{\mu}=  1+\frac{2s}{-\tau_-}.
 \end{equation}

 Let $\tau_0=\tau_-$,
 which verifies that for $x\in \R^N\setminus  \bar B_r$
$$
\cL^s_\mu  u_0(x)  \geq      d_0^p |x|^{p\tau_-+\theta}= d_0^p |x|^{\tau_1-2s}\quad {\rm in}\ \, \R^N\setminus  \bar B_r,
$$
where
$$\tau_1:=p\tau_0+ 2s.$$

If $p\tau_0 \geq   \tau_+-2s$, then
$$  u_0^p\geq  d_0^p  |x|^{\tau_+-2s}$$
and a contradiction follows by Lemma \ref{re 2.1}. We are done.  \smallskip

If not, by Lemma \ref{teo 2.2}, we have that
$$u_0(x)\geq d_1|x|^{\tau_1}\quad {\rm in}\ \, \R^N\setminus   \bar B_r.$$
Iteratively, we recall that
$$\tau_j:=p\tau_{j-1}+\theta+2s,\quad j=1,2,\cdots .$$
Note that for $p\in(0,   p^*_{\theta,\mu})$
$$\tau_1-\tau_0=(p-1)\tau_0+\theta+2s >0.$$
 If $\tau_{j+1}=\tau_jp+\theta+2s\in  (\tau_-,\tau_+)$,
it following by Theorem \ref{teo 2.2}  that
$$u_0(x)\geq d_{j+1} |x|^{\tau_{j+1}},$$
where
$$\tau_{j+1}=p\tau_j+2s+\theta>\tau_j.$$
  If $p\tau_{j+1}+\theta\geq  \tau_-$, we are done by lemma \ref{teo 2.2}.  In fact, this iteration could stop by finite times since  $\tau_j\to+\infty$ as $j\to+\infty$ if $p\geq 1$.\smallskip

  {\bf Step 2.}  Nonexistence in the critical case $\mu>\mu_0$  $p=p^*_{\theta,\mu}>1$.
  Note  that for some $\sigma_0>0$
 $$\frac12  u_0^{p-1}(x)\geq \sigma_0|x|^{-2s} \quad{\rm in}\ \, \R^N\setminus B_r.$$
 Here we can assume that $\mu-\sigma_0\geq \mu_0$, otherwise, we only take a smaller value for $\sigma_0$.
 So we can write problem (\ref{eq 1.1}) as following
\begin{equation}\label{eq 1.1 new-ex}
 \cL_{\mu-\sigma_0}^s u_0\geq \frac12    u_0^p \quad {\rm in}\ \  R^N\setminus B_r,
 \end{equation}
 which the critical exponent
 $$p^*_{\mu-\sigma_0 }=  1+\frac{2s }{-\tau_-(s,\mu-\sigma_0)} >    1+\frac{2s }{-\tau_+(s,\mu)}=p^*_{\mu}, $$
 where $\mu\in(\mu_0,+\infty)\mapsto \tau_+(s,\mu)$ is strictly increasing.
Thus a contradiction comes from {\bf step 1 } for (\ref{eq 1.1 new-ex}).
  \hfill$\Box$  \medskip

\noindent {\bf Proof of Proposition \ref{teo 3}. } It is the particular case $\mu=0$ and $p=p^*$ in Theorem \ref{teo 3}.  \hfill$\Box$

    \bigskip\bigskip

   \noindent{\bf \small Acknowledgements:} {\footnotesize This work is supported by NNSF of China, No: 12071189, by the
Jiangxi Provincial Natural Science Foundation, No: 20212ACB211005 and 20202ACBL201001.}


\begin{thebibliography}{99}
%\addtolength{\itemsep}{-0.5ex}



  \bibitem {AC}  W. Ao, H. Chan, A. Dela-Torre,  M. Fontelos,   G. Mar\'ia del Mar, J. Wei.  On higher-dimensional singularities for the fractional Yamabe problem: a nonlocal Mazzeo-Pacard program. {\it Duke Math. J. 168(17)},  3297--3411 (2019).

   \bibitem {AD} W. Ao,  A.  DelaTorre,  M. Gonz\'alez,  J. Wei.  A gluing approach for the fractional Yamabe problem with isolated singularities. {\it J. Reine Angew. Math. 763},  25--78 (2020).

  \bibitem {A1} P. Aviles. On isolated singularities in some nonlinear partial differential equations.
  {\it Indiana Univ. Math. J. 35}, 773--791 (1983).

    \bibitem {A2} P. Aviles. Local behavior of solutions of some elliptic equations. {\it Commum. Math. Phys. 108},  177--192 (1987).


  \bibitem {BV}  M.-F. Bidaut-V\'eron, L. V\'eron. Nonlinear elliptic equations on compact Riemannian manifolds and asymptotics of Emden equations. {\it Invent. Math. 106 (3)}
489--539 (1991).


  \bibitem {BM} K. Bogdan, B. Dyda. The best constant in a fractional Hardy inequality.
{\it  Math. Nachr. 284 (5-6),} 629--638 (2011).

 \bibitem {CGS} L. Caffarelli, B. Gidas, J. Spruck. Asymptotic symmetry and local behavior of semilinear elliptic equations with critical Sobolev growth. {\it Comm. Pure Appl. Math. 42} 271--297 (1989).

 \bibitem {CJ}  L. Caffarelli, T. Jin, Y. Sire, J. Xiong. Local analysis of solutions of fractional semi-linear elliptic equations with isolated singularities. {\it Arch. Ration. Mech.
Anal. 213(1)},  245--268 (2014).

  \bibitem {CS0} L. Caffarelli,  L. Silvestre.  Regularity theory for fully nonlinear integro-differential equations. {\it Comm. Pure Appl. Math. 62(5)},  597--638 (2009).

  \bibitem {CS1} L. Caffarelli,  L. Silvestre. Regularity results for nonlocal
equations by approximation. {\it Arch. Ration. Mech. Anal.  200(1),}
59--88 (2011).


\bibitem{CS2} L. Caffarelli,  L. Silvestre. An extension problem related to the fractional Laplacian.  {\it Comm. Part. Diff. Eq. 32}, 1245--1260 (2007).

%{\color{blue} \bibitem {CD} H. Chan, A. DelaTorre. From fractional Lane-Emden-Serrin equation -- existence, multiplicity and local behaviors via classical ODE -- to fractional Yamabe metrics with singularity of "maximal" dimension, {\it arXiv:2109.05647.} (2021)}


\bibitem {CFQ} H. Chen,  P. Felmer, A. Quaas. Large solutions to elliptic equations involving fractional Laplacian. {\it Ann. Inst. Henr.  Poincar\'e-AN, 32,} 1199--1228 (2015).

\bibitem {CQ} H. Chen,   A. Quaas.   Classification of isolated singularities of nonnegative solutions to fractional semilinear elliptic equations and the existence results. {\it J. London Math. Soc.   97(2),} 196--221 (2018).

 \bibitem {CV1} H. Chen,  L. V\'{e}ron, Semilinear fractional elliptic equations involving measures. {\it J. Diff. Eq. 257(5)}, 1457--1486 (2014).

\bibitem {CW} H. Chen, T. Weth. The Poisson problem for the fractional Hardy operator: Distributional identities and singular solutions.  {\it Trans.  Amer. Math. Soc. 374(10)},   6881--6925   (2021).


  \bibitem{CLO}  W. Chen,  C. Li,  B. Ou.  Qualitative properties of solutions for an integral equation. {\it Discr. Cont. Dyn. Syst. 12(2),} 347--354 (2005).

 \bibitem{CLO1}  W. Chen,  C. Li,  B. Ou.  Classification of solutions for an integral equation. {\it Commun. Pure Appl. Math. 59(3),} 330--343 (2006).
 
  \bibitem {CS} Z. Chen,   R. Song. Estimates on Green functions and poisson kernels for symmetric stable process. {\it Math.
Ann. 312}, 465--501 (1998). 

 \bibitem{DdG}   A. DelaTorre, M. del Pino, M. Gonz\'alez, J. Wei.
 Delaunay-type singular solutions for the fractional Yamabe problem. {\it Math. Ann. 369}, 597--626 (2017) .

 \bibitem{DG} A. DelaTorre,  M. Gonz\'alez.  Isolated singularities for a semilinear equation for the fractional Laplacian arising in conformal geometry. {\it Rev. Mat. Iberoam. 34(4)}, 1645--1678 (2018).

  \bibitem{NPV}  E. Di Nezza, G. Palatucci, E. Valdinoci.  Hitchhiker's guide to the fractional Sobolev spaces.
  {\it Bull. Sci. Math. 136}, 521--573 (2012).

  \bibitem {FW}   P. Felmer,  Y. Wang.  Radial symmetry of positive solutions to equations involving the fractional Laplacian.
 {\it Comm. Cont. Math. 16 (1),}   1350023, 24 pp  (2014).

 \bibitem {FQ} P. Felmer, A. Quaas. Fundamental solutions and Liouville type theorems for nonlinear integral operators. {\it Adv.  Math. 226}, 2712--2738 (2011).


  \bibitem {GS} B. Gidas, J. Spruck. Global and local behavior of positive solutions of nonlinear elliptic equations. {\it Comm. Pure Appl. Math. 34}, 525--589 (1981).


  \bibitem {HLT}  Z.  Han, Y.  Li,  E. Teixeira. Asymptotic behavior of solutions to the k-Yamabe equation near isolated singularities. {\it Invent. Math.  182},   635--684 (2010).


 \bibitem{JLX} T. Jin,  Y. Li,  J. Xiong.   On a fractional Nirenberg problem, part I: blow up analysis and compactness of solutions. {\it J. Eur. Math. Soc.  16(6),}   1111--1171 (2014).

\bibitem{KM} N. Korevaar, R. Mazzeo, F. Pacard, R. Schoen. Refined asymptotics for constant scalar curvature metrics with isolated singularities. {\it Invent. Math. 135},  233--272 (1999).

% \bibitem{JX} T. Jin,  J. Xiong,   Asymptotic symmetry and local behavior of solutions of higher order conformally invariant equations with isolated singularities,  {\it Ann. Inst. H. Poincar\'e--AN 38(4)},  1167--1216 (2021).

 \bibitem{L}  P.-L. Lions. Isolated singularities in semilinear problems. {\it J. Diff. Eq.
38},  441--450 (1980).



% \bibitem {R}  W. Rudin, Function analysis, {\it McGraw-Hill} (1973).

 %\bibitem {S} L. Schwartz, Theorie des distributions, {\it Hermann, Paris} (1966).

  \bibitem{P}  F. Pacard. Existence and convergence of positive week solutions of
  $-\Delta u=u^{\frac{n}{n-2}}$ in bounded domain of $\R^n$, $n\geq3$.  {\it Calc. Var. 1}, 243--265 (1993).


\bibitem {RS} X. Ros-Oton, J. Serra. The Dirichlet problem for the fractional
laplacian: regularity up to the boundary. {\it J. Math. Pures Appl. 101},  275--302 (2014).

  \bibitem{S0} R. Schoen. The existence of weak solutions with prescribed singular behavior for a conformally invariant scalar equation. {\it Comm. Pure Appl. Math., 41},  317--392 (1988).




 \bibitem{S}  J, Serrin. Local behavior of solutions of quasi-linear equation. {\it Acta. Math. 111}, 247--302 (1964).
 
  \bibitem{SV}  R. Servadei,  E. Valdinoci.  Variational methods for non-local operators of elliptic type. {\it Disc. Contin. Dyn. Syst. 33(5)},   2105--2137 (2013). 

 \bibitem {Ls07} L. Silvestre.  Regularity of the obstacle problem for a fractional power of the Laplace operator.  {\it Comm. Pure   Appl. Math. 60 (1),} 67--112 (2007).


 \bibitem{TX} J. Tan, J. Xiong. A Harnack inequality for fractional Laplace equations with lower order terms.  {\it Dis. Cont. Dyn. Syst. 31(3)}, 975--983 (2011).

  \bibitem {WW}   J. Wei, K. Wu. Local behaviour of solutions to a fractional equation with isolated behavior of solutions to a fractional equation with isolated singularity and critical Serrin exponent. {\it Disc. Cont. Dyn. Syst.   42(8),}  4031--4050 (2022).  

 \bibitem{YZ} H. Yang, W. Zou. On isolated singularities of fractional semilinear elliptic equations.
 {\it  Ann.  Inst. Henri Poincar\'e--AN 38(2),}  403--420  (2021).

% \bibitem{V} L. V\'eron.   Singularities of solutions of second order quasilinear equations,  {\it Pitman Research Notes in Mathematics Series, 353. Longman, Harlow}  viii+377 pp. (1996).


\end{thebibliography}
 \end{document}